\documentclass[draft,onecolumn]{IEEEtran}
\usepackage[latin9]{inputenc}
\usepackage{geometry}
\usepackage{color}
\usepackage{float}
\usepackage{amsmath}
\usepackage{amsthm}
\usepackage{amssymb,bbm}
\usepackage{graphicx}

\usepackage{flushend}

\usepackage{mathtools}
\mathtoolsset{showonlyrefs=true}

\usepackage[linesnumbered,ruled,vlined,noend]{algorithm2e}

\makeatletter
\usepackage{color}

\DeclareMathOperator*{\argmax}{arg\,max}
\DeclareMathOperator*{\argmin}{arg\,min}

\newcommand{\bone}{\boldsymbol{1}}
\newcommand{\round}[1]{[#1]_r}

\usepackage{babel}

\newtheorem{lemma}{Lemma}
\newtheorem{problem}{Problem}
\newtheorem{proposition}{Proposition}
\newtheorem{theorem}{Theorem}
\newtheorem{definition}{Definition}
\newtheorem{assumption}{Assumption}
\newtheorem{remark}{Remark}
\newtheorem{example}{Example}

\begin{document}
\title{Distributed Asynchronous Mixed-Integer Linear Programming with Feasibility Guarantees}

\author{Luke Fina$^{1}$, Christopher Petersen$^{1}$, and Matthew Hale$^{2}$\thanks{$^{1}$ Department of Mechanical and Aerospace Engineering, University
of Florida, Gainesville, FL USA. Emails: \texttt{\{l.fina,c.petersen1\}@ufl.edu}. \\
$^{2}$ School of Electrical and Computer Engineering,
Georgia Institute of Technology, Atlanta, GA USA.
Email: \texttt{matthale@gatech.edu}. \\
This work was supported by AFOSR under grant FA9550-19-1-0169, ONR
under grants N00014-19-1-2543, N00014-21-1-2495, and N00014-22-1-2435,
and by AFRL under grant
FA8651-08-D-0108 TO48.}}

\maketitle

\begin{abstract}
In this paper we solve mixed-integer linear programs (MILPs) via distributed
asynchronous saddle point computation. To solve a MILP, we relax it with a linear program approximation. 
We first show that if the linear program relaxation satisfies Slater's condition, then 
relaxing the problem, solving it, and rounding the relaxed solution 
produces a point that is guaranteed to satisfy the constraints of the original MILP. 
Next, we form a Lagrangian saddle point problem
that is equivalent to the linear program relaxation, and 
then we regularize the Lagrangian in both the
primal and dual spaces. We then develop a parallelized algorithm
to compute saddle points of the regularized Lagrangian, and 
we show that it is tolerant to 
asynchrony in the computations and communications of primal and dual
variables. Suboptimality bounds and convergence rates are presented
for convergence to a saddle point. Simulation results illustrate these theoretical
developments in practice, and show that relaxation and regularization combined
typically have only a mild impact on the suboptimality of the solution obtained.     
\end{abstract}

\section{Introduction\label{sec:Introduction}}

Numerous problems in autonomy can be formulated as mixed-integer
programs in which the integer constraints lead to a more realistic
model of the system than a problem without integer constraints. In
particular, mixed-integer programs have been used to model piece-wise
affine functions which can approximate nonlinear dynamical systems
for control~\cite{richards2005mixed}, deep neural networks to improve robustness~\cite{fischetti2018deep},
task assignment problems~\cite{hendricksondecentralized}, and trajectory
planning problems with collision avoidance constraints~\cite{ioan2021mixed}.
Mixed-integer programs are useful in their
modeling power, but scale poorly in comparison to continuous optimization
problems due to the NP-Hardness of solving problems with integrality
constraints~\cite{conforti2014integer}. 

Therefore, methods have been
developed to computationally reach approximate solutions. 
Many existing
methods are centralized, but as the size of problems increases, 
decentralized schemes may be desirable to partition problems into
smaller pieces and accelerate computations. To help address this need, in this work we solve mixed-integer
linear programs (MILPs) using a network of agents that 
is permitted to compute and communicate
information asynchronously. 

Our approach relies in part on relaxing MILPs, which is done by
removing integrality constraints, solving the 
resulting (ordinary) linear program (LP), and rounding
its solution. 
Existing work shows that this approach produces a point that
satisfies the constraints of the original MILP if that MILP
is ``granular''~\cite{neumann2019feasible}. 
Granularity of a MILP is a certain
technical condition that can be difficult to 
check in general. 
However, we circumvent this issue by providing a new result that
shows that if the LP relaxation 
satisfies Slater's condition, 
then the granularity condition is satisfied, and
therefore rounding a solution to the relaxed LP 
is guaranteed to produce a point that satisfies the constraints of the original MILP.
This result eliminates the need to directly verify granularity, and it also guarantees 
feasibility, which has been difficult to guarantee
in some prior works.

For relaxed problems that satisfy Slater's condition, 
we find a feasible approximate solution that has bounded,
quantifiable suboptimality. 
 We do this by translating the relaxed LP to an equivalent
Lagrangian saddle point problem, and we regularize this Lagrangian 
in the primal
and dual spaces to make it strongly convex-strongly concave.
Regularization ensures uniqueness of a solution, i.e., it removes degeneracy by construction.
Regularization can also perturb solutions, and we develop conditions under which
feasibility of a rounded solution is still ensured under regularization. 

We also show that regularization provides robustness to asynchrony.
Asynchrony is inherent to some settings, such as large networks
of many processors and/or multi-agent settings with limited and congested
bandwidth. To address these cases, we develop
an algorithmic framework that is robust to asynchrony in communications
and computations. 
The algorithm we develop is a parallelized saddle point finding algorithm
for the regularized Lagrangian. This algorithm
partitions problems into blocks that are either scalars or sub-vectors
of the primal or dual decision variables, and this algorithm is shown
to be tolerant to bounded delays in computations and communications. The
block structure reduces the required
communications among agents
because linearity implies that the primal variables are separable, which implies
that agents can compute their local gradients without needing to share
information with each other about their primal variables. 
The required communications depend upon the coupling
induced by the constraints in a problem in a way that we make precise. 
We present convergence analysis for our
algorithm, along with a suboptimality bound between the solution to
the original MILP and the approximate solution obtained using our
algorithm. 

Our developments draw on recent advances in distributed
primal-dual optimization and solving granular mixed-integer linear
programs. Similar work appears in \cite{camisa2021distributed,camisa2019primal},
but one major distinction is that those works use a "large-scale'' mixed-integer linear program formulation 
in which the number of constraints is small relative to the number
of decision variables. 
The problem formulation that we consider is more general because it does not rely on problems having a small number of coupling constraints, and it is computationally simpler since our method
does not require the construction of a convex hull for the mixed-integer constraint set, which
is nontrivial in general~\cite{conforti2014integer}. 

Distributed primal-dual optimization is a foundational
element for the current work because it is used to solve the relaxed
MILP \cite{hendrickson2021totally}. Earlier works on decentralized
primal-dual constrained optimization address the effect of delay bounds,
but \cite{hale2017asynchronous} did not address integer programming
and \cite{falsone2017dual} used consensus-based methods while we use a block-based approach. 
In multi-agent
optimization, several recent works on MILPs
draw from~\cite{vujanic2016decomposition}, specifically \cite{falsone2017dual,falsone2020distributed,camisa2019primalasync}. Our work focuses on a more general class of MILPs than \cite{vujanic2016decomposition} and subsequent derived works.
The asynchronous solution of MILPs was also addressed in \cite{bragin2020distributedasync}, but our work obtains a feasible mixed-integer solution without decomposing the original problem into subproblems.

To summarize, the main contributions of this paper are the following:
\begin{itemize}
\item We show that satisfaction of Slater's condition guarantees that the rounded solution to a LP relaxation 
will satisfy the constraints of the original MILP (Theorems \ref{thm:SlaterConnectionToGranularity} and~\ref{thm:saddle_feasible}).
\item We develop a distributed asynchronous algorithm for solving MILPs that guarantees constraint satsifaction of the final answer obtained (Algorithm \ref{alg:algorithm1}).
\item We bound the suboptimality gap between the original
MILP and a rounded solution
to its corresponding relaxed and 
regularized Lagrangian saddle point problem (Theorem
\ref{thm:suboptimalityEB}).
\item We present a rate of convergence for our algorithm under asynchrony
(Theorems \ref{thm:dualConvergence} and \ref{thm:overallConvergence}).
\end{itemize}

The rest of this paper is organized as follows. Preliminaries and problem statements
 are given in Section \ref{sec:Preliminaries-and-Problem}. 
 Our main algorithm is given in Section \ref{sec:algprelims}. Section
 \ref{sec:Distributed-Algorithm} presents convergence
 analysis and suboptimality bounds. Section \ref{sec:Simulation}
 demonstrates the algorithm in simulation, and Section \ref{sec:Conclusion}
 concludes. \\

\noindent\textbf{Notation and Terminology} We write ``LP'' for ``linear program'' and ``MILP'' for ``mixed-integer linear program.'' We use~$\mathbb{R}$ to denote the real numbers,
$\mathbb{R}_{+}$ to denote the non-negative real numbers, 
$\mathbb{N}$ to denote the natural numbers, 
$\mathbb{N}_0$ to denote the natural numbers including zero, and~$\mathbb{Z}$ for the integers. 
For~$n \in \mathbb{N}$ we use~$[n] = \{1, \ldots, n\}$, and
we also use~$a\mathbb{Z} = \{az \mid z \in \mathbb{Z}\}$ for
any~$a \in \mathbb{R}$. 
For a vector~$v \in \mathbb{R}^n$, we write~$v^T$ for its transpose,
and we use~$\|\cdot\|$ (without a subscript) to denote
the Euclidean norm. 
We will interchangeably use~$\mathbb{R}^n \times \mathbb{R}^m$ and~$\mathbb{R}^{n+m}$, and
similar for other trivially isomorphic sets. 
The symbol~$\boldsymbol{1}$ denotes the vector of all~$1$'s. 
We use~$\lfloor \cdot \rfloor$ to denote the floor function and~$\round{\cdot}$ to denote rounding to the nearest integer;
ties between two closest integers can be broken arbitrarily. 

\section{Formulating Relaxed MILPs}  \label{sec:Preliminaries-and-Problem}
We solve the following form of MILPs. 
\begin{problem}[Original MILP] \label{prob:Original-MILP}
    
    \begin{align*}
    \underset{(x,\tilde{y})\in X \times\mathbb{Z}^{m}}{\textnormal{minimize}} & \enskip a^{T}x+d^{T}\tilde{y}\\
    \textnormal{subject to} & \enskip Ex+F\tilde{y}\leq h, \\
     & y^{\ell}\leq \tilde{y}\leq y^{u}, 
    \end{align*}
    where $X\subset \mathbb{R}^{n}$ is a polyhedral set, $a\in \mathbb{R}^{n}, d \in \mathbb{R}^{m}, E\in \mathbb{R}^{p\times n}, F \in \mathbb{Z}^{p\times m}, h\in\mathbb{R}^{p}$, $y^{\ell}\in\mathbb{Z}^{m}$ is a lower bound  on $\tilde{y}$, and $y^{u}\in\mathbb{Z}^{m}$ is an upper bound  on $\tilde{y}$.
\end{problem}

    We enforce the following assumption about~$X$. 

    \begin{assumption} \label{as:Xcompact}
        The set~$X \subseteq \mathbb{R}^n$ is non-empty, compact, and convex. It can be decomposed into
        the Cartesian product~$X = X_1 \times \cdots X_n$. 
    \end{assumption}

   Assumption~\ref{as:Xcompact} ensures that the feasible
    region of Problem~\ref{prob:Original-MILP} is bounded, and it will be used 
    to develop a projected update law in Algorithm \ref{alg:algorithm1}.

    We will solve problems of the form of Problem~\ref{prob:Original-MILP} by relaxing the integer part, i.e., by optimizing over~$y \in \mathbb{R}^m$ instead of~$\tilde{y} \in \mathbb{Z}^m$,
    then solving this relaxed problem, and rounding its solution.
    It is well-known, e.g.,~\cite{wolsey2020integer}, 
    that in general this procedure does not produce points that satisfy the constraints
    in the original MILP. 
    However, for the class of mixed-integer linear programs called ``granular'' mixed-integer programs, 
    solving a continuous relaxation and rounding its solution is
    guaranteed to produce a point that satisfies the constraints of the original MILP~\cite{neumann2019feasible}. 
    
    Defining granularity first requires
    a few definitions from~\cite{neumann2019feasible}. 
    We define the feasible region in Problem~\ref{prob:Original-MILP} as 
    \begin{equation} \label{eq:mmilp}
    M_{MILP}=\left\{ (x,\tilde{y})\in X\times\mathbb{Z}^{m} :  Ex + F\tilde{y} \leq h,  y^{\ell} \leq \tilde{y} \leq y^u \right\},
    \end{equation}    
    where $E$, $F$, $y^{\ell}$, $y^u$, and $h$ are from Problem \ref{prob:Original-MILP}. 
    We use~$f_i^T$ to denote the~$i^{th}$ row of~$F$
    and~$e_i^T$ to denote the~$i^{th}$ row of~$E$ 
    for~$i \in \{1, \ldots, p\}$. For each $i\in \{1,...,p\}$,
    we define $\omega_{i}$ as follows.
    If~$e_i \neq 0$, then~$\omega_i = 0$.
    If~$e_i = 0$, then~$\omega_i$ is set equal to 
    the greatest common divisor of the non-zero entries of $f_{i}\in \mathbb{Z}^{m}$. Then we set
    \begin{equation}\label{eq:omegaVector}
    \omega := \big(\omega_{1},\ldots,\omega_{p}\big)^{T}\in \mathbb{N}^{p}_{0}.
    \end{equation}
    Using~$h_i$ to denote the~$i^{th}$ entry of~$h$
    for each~$i \in \{1, \ldots, p\}$, we further define 
    \begin{equation}\label{eq:hOmega}
        \lfloor h_{i} \rfloor_{\omega_{i}}  = \max \{ t\in \omega_{i}\mathbb{Z} \mid t\leq h_{i}\},     
    \end{equation}
    where we set~$\lfloor h_i \rfloor_{0} = h_i$.
    Then, for~$h \in \mathbb{R}^p$ we set  
    \begin{equation} \label{eq:hOmegaFull}
        \lfloor h \rfloor_{\omega} := \big(\lfloor h_{1} \rfloor _{\omega_{1}} ...  \lfloor h_{p} \rfloor _{\omega_{p}}\big)^{T}.
    \end{equation}
    
    Using these terms, the mixed-integer linear feasible set will be constructed using a larger auxiliary set.
    
    \begin{definition}[Section~2.2 of~\cite{neumann2019feasible}] \label{def:EIPS}
    The ``enlarged inner parallel set'' associated with Problem~\ref{prob:Original-MILP} is 
    \begin{equation} \label{eq:EIPS}
    M_{\xi}=\left\{ (x,y)\in X\times\mathbb{R}^{m} \mid Ex + Fy \leq \lfloor h\rfloor_{\omega}+\xi\omega-\frac{1}{2}\rho, y^{\ell} + \left(\frac{1}{2} - \xi\right)\bone \leq y \leq 
    y^u + \left(\xi - \frac{1}{2}\right)\bone
    \right\},
    \end{equation}
    where~$\xi \in (0,1)$, 
    $\lfloor h \rfloor_{\omega}$ is from \eqref{eq:hOmegaFull}, $\omega$ is from \eqref{eq:omegaVector}, $\rho := [\|f_{1} \|_{1}, ...,\|f_{_{p}} \|_{1}]^{T}$ where~$f_i^T$ is the~$i^{th}$ row of~$F$, 
    and where the inequalities in~\eqref{eq:EIPS} are to be interpreted
    componentwise. 
    \end{definition}

    It is clear that for~$\xi_1, \xi_2 \in (0,1)$ we have
    \begin{equation} \label{eq:Mcontainment} 
        \xi_1 \leq \xi_2 \quad\textnormal{ implies }\quad M_{\xi_1} \subseteq M_{\xi_2}.
    \end{equation}

    The enlarged inner parallel set can be used as the feasible region
    for a relaxed MILP that considers~$y \in \mathbb{R}^m$ 
    rather than~$\tilde{y} \in \mathbb{Z}^m$,
    and, under mild conditions, rounding the solution
    to that relaxed MILP produces a point that is feasible
    for the original MILP. 
    
    \begin{proposition}[{\cite[Theorem 2.6]{neumann2019feasible}}] \label{prop:GranRounding}
        For any $\xi \in [\xi_{e},1)$, any rounding $(x,\round{y})$ of any $(x,y)\in M_{\xi}$ lies in $M_{\text{MILP}}$,
     where 
    \begin{equation} \label{eq:xiedef}
        \xi_{e} = \begin{cases} 
        \max\left\{ \frac{h_{i}-\lfloor h_{i}\rfloor_{\omega_{i}}}{\omega_{i}}\mid i=1,...,p,\omega_{i}\neq0\right\} & \omega\neq0\\
        0 & \textnormal{else}
        \end{cases}
        \end{equation}
        satisfies $\xi_e < 1$. 
    \end{proposition}
    

    
    Under Proposition \ref{prop:GranRounding}, we next construct a linear program 
    whose solution can be rounded to produce a point that satisfies the constraints
    of Problem~\ref{prob:Original-MILP}. 
    %
   

    \begin{problem}[Relaxed LP] \label{prob:matrix_Feasible_LP}
        Fix~$\xi \in [\xi_e, 1)$ for~$\xi_e$ as defined in~\eqref{eq:xiedef}. Then
        \begin{align*}
            \underset{(x,y) \in X \times \mathbb{R}^m}{\text{minimize}} 
            & \enskip a^{T}x + d^Ty \\
            \text{subject to} 
            & \enskip g_{\xi}(x, y) \!=\! 
            \left(\!\!\!\begin{array}{c} 
             Ex + Fy - \lfloor h\rfloor_{\omega} - \xi\omega + \frac{1}{2}\rho \\
             y - y^u - \left(\xi - \frac{1}{2}\right)\bone \\
             -y + y^{\ell} + \left(\frac{1}{2} - \xi\right)\bone
            \end{array}\!\!\!\right) \!\leq\! 0,
            \end{align*}
    where 
    $a$,~$d$,~$X$~$E$,~$F$,~$y^{\ell}$, and~$y^u$ 
    are from Problem~\ref{prob:Original-MILP}, 
    $\lfloor h \rfloor_{\omega} \in\mathbb{R}^{p}$ 
    is from \eqref{eq:hOmegaFull}, 
    $\omega$ is from \eqref{eq:omegaVector}, and $\rho := [\|f_{1} \|_{1}, ...,\|f_{_{p}} \|_{1}]^{T}$,
    where~$f_i^T$ is the~$i^{th}$ row of~$F$. 
    \end{problem}

    We note that the function~$g_{\xi}$ in Problem~\ref{prob:matrix_Feasible_LP}
    encodes the same feasible region as the 
    enlarged inner parallel
    set~$M_{\xi}$ in~\eqref{eq:EIPS}. The existence of a solution to Problem~\ref{prob:matrix_Feasible_LP} is guaranteed if the relaxed feasible set~$M_{\xi}$ is nonempty for the choice of $\xi$ that is used to formulate it. In some problems, 
    the enlarged inner parallel set~$M_{\xi}$ could be empty for 
    all values of~$\xi$, even for 
    a non-empty set $M_{MILP}$, and the following definition is used to exclude such cases.
    
    \begin{definition}[Definition~2.7 in~\cite{neumann2019feasible}] \label{def:MILPGran}
        The feasible set $M_{\text{MILP}}$ in~\eqref{eq:mmilp} is ``granular'' if its enlarged inner parallel set $M_{\xi}$ in~\eqref{eq:EIPS} 
        is nonempty for some $\xi \in [\xi_{e},1)$. Moreover, Problem \ref{prob:Original-MILP} is said to be
        ``granular'' if its feasible set $M_{\text{MILP}}$ is \textit{granular}.
    \end{definition}
    
    We also assume the following. 
    %
    %
    %
    %

    \begin{assumption}[Slater's Condition for~$g_{\xi_e}$] \label{as:slater}
        There exists~$(\bar{x}, \bar{y}) \in X \times \mathbb{R}^m$ such that~$g_{\xi_e}(\bar{x}, \bar{y}) < 0$, where~$\xi_e$ is from 
        Proposition~\ref{prop:GranRounding}. 
    \end{assumption}
     
     Assumption~\ref{as:slater} will be used to ensure that strong duality
     holds for the problems that we consider.

    \begin{remark} \label{rem:slater}
        Due to~\eqref{eq:Mcontainment}, the satisfaction
        of Assumption~\ref{as:slater} by~$g_{\xi_e}$ also
        implies that~$g_{\xi}(\bar{x},\bar{y}) < 0$ for
        all~$\xi \in [\xi_e, 1)$. That is, 
        satisfaction of Assumption~\ref{as:slater} by~$g_{\xi_e}$
        implies that Slater's condition holds
        for all~$g_{\xi}$ with~$\xi \in [\xi_e, 1)$. 
    \end{remark}

    Next, we show that 
    under Assumption~\ref{as:slater}, 
    one can solve Problem~\ref{prob:matrix_Feasible_LP}
    for any value of~$\xi \in [\xi_e, 1)$ and round its solution to
    obtain a point that is feasible for Problem~\ref{prob:Original-MILP}.

\begin{theorem}[Slater's Condition Implies Granularity]  \label{thm:SlaterConnectionToGranularity}
    Consider Problem~\ref{prob:Original-MILP} and its relaxation to Problem~\ref{prob:matrix_Feasible_LP}, and 
    suppose that Assumptions~\ref{as:Xcompact} and~\ref{as:slater} are
    satisfied by Problem~\ref{prob:matrix_Feasible_LP}. 
    Then the original mixed-integer linear program in Problem~\ref{prob:Original-MILP} is granular. Furthermore, 
    at least one solution~$(x^*, y^*)$ exists for
    Problem~\ref{prob:matrix_Feasible_LP}, and, for
    every such solution, $(x^*,\round{y^*})$ 
    satisfies the constraints in Problem~\ref{prob:Original-MILP}.
\end{theorem}

\emph{Proof: } See Appendix~\ref{ss:theorem1}. \hfill $\blacksquare$

Past work \cite{glover1975pitfalls} has shown that it can be difficult
in general to ensure that rounded solutions to relaxed integer programs 
satisfy the constraints in the original problem. 
Methods have been developed to provide such assurances \cite{bertsimas2005optimization}, 
which include various forms of heuristics. Alternatively,
Theorem~\ref{thm:SlaterConnectionToGranularity} provides a guarantee
that rounded solutions to relaxed problems are feasible with respect
to the constraints of the original MILP, and no heuristics or post-processing
of such solutions are required beyond rounding. 

Slater's condition is commonly used in the nonlinear programming literature,
and one can verify that it holds, e.g., by solving a simple
feasibility program \cite{bertsekas2003convex}. Here, 
one can use it to ensure that relaxing and rounding for MILPs will
provide a feasible solution. While existing work has used
Slater's condition to provide error bounds for rounded
solutions to relaxed problems, we are not aware of prior work
that connects it to the property of granularity.

\section{Lagrangian Setup and Analysis} \label{sec:algprelims}
In this section we reformulate Problem~\ref{prob:matrix_Feasible_LP} into
an equivalent Lagrangian saddle point problem. Then we characterize
the relationship between solutions to such saddle point problems
and solutions to the original MILPs. Finally,
we regularize these Lagrangians to 
make them robust to being solved
asynchronously in Section~\ref{sec:Distributed-Algorithm}. 

\subsection{Lagrangian Saddle Point Formulation}

    
We will solve Problem \ref{prob:matrix_Feasible_LP} using a primal-dual approach.
Specifically, we will use a distributed asynchronous form of the classic Uzawa
iteration~\cite{arrow1958studies}, 
also called gradient descent-ascent, to find a saddle point
of the Lagrangian associated with Problem \ref{prob:matrix_Feasible_LP}.
In this section, we assume that some~$\xi \in [\xi_e, 1)$ has
been selected and remains fixed. 
For concise notation in the forthcoming
saddle point problem, we introduce several new symbols.
Namely, we set
\begin{align}
c &= (a^T, d^T)^T \in \mathbb{R}^{n+m} \label{eq:c} \\
A &= (E \,\,\,\, F) \in \mathbb{R}^{p \times (n + m)} \label{eq:A} \\
b &= \lfloor h \rfloor_{\omega} \label{eq:b} \\
\nu &= \xi\omega - \frac{1}{2}\rho \label{eq:nu} \\
z &= (x^T, y^T)^T \in \mathbb{R}^{n + m} \label{eq:z} \\
Z &= X \times Y_{\xi} \label{eq:Z} \\
\bar{z} &= (\bar{x}^T, \bar{y}^T)^T, \label{eq:zbar} 
\end{align}
where~$a$, $d$, $E$, $F$, $h$, $\omega$, $\xi$, $\rho$, $x$, $y$, and~$X$
are from Problem~\ref{prob:matrix_Feasible_LP}, 
the points~$\bar{x}$ and~$\bar{y}$ are from Assumption~\ref{as:slater},
and
where
\begin{equation}
Y_{\xi} = \left\{y \in \mathbb{R}^m \mid 
        y^{\ell} + \left(\frac{1}{2} - \xi\right)\bone \leq y \leq 
    y^u + \left(\xi - \frac{1}{2}\right)\bone
        \right\}.
\end{equation}

\begin{remark} \label{rem:yprod} 
The set~$Y_{\xi}$ admits a Cartesian product structure without any additional assumptions.
Indeed, we see that~$y_i$ must satisfy~$y^{\ell}_i + (\frac{1}{2} - \xi) \leq y_i \leq y^u_i + (\xi - \frac{1}{2})$.
Thus, by defining the interval
\begin{equation}
Y_i = \left[y^{\ell}_i + (\frac{1}{2} - \xi), y^u_i + (\xi - \frac{1}{2})\right] \subseteq \mathbb{R},
\end{equation}
we see that~$y_i \in Y_i$ for all~$i \in [m]$ if and only if~$y \in Y_{\xi}$. 
\end{remark}

The Lagrangian 
associated with Problem~\ref{prob:matrix_Feasible_LP} 
is 
\begin{equation}
L(z,\lambda)=c^{T}z+\lambda^{T}(Az-b-\nu), 
\end{equation}
where~$z \in Z$ and where~$\lambda \in \mathbb{R}_{+}^p$ is a vector of Karush-Kuhn-Tucker (KKT)
multipliers. This~$L$
is affine in both $z$ and $\lambda$, and thus it is convex in
$z$ and concave in $\lambda.$ To ensure convergence of a decentralized
asynchronous implementation of the Uzawa iteration, we will regularize
$L$ to make it strongly convex in $z$ and strongly concave in $\lambda$.
A Tikhonov regularization adds a quadratic
term in~$z$ and subtracts a quadratic
term in~$\lambda$ to form the regularized
Lagrangian
\begin{equation} \label{eq:Lregdef}
L_{\alpha,\delta}(z, \lambda) = L(z, \lambda) + \frac{\alpha}{2}\|z\|^2 - \frac{\delta}{2}\|\lambda\|^2,
\end{equation}
where~$\alpha > 0$ and~$\delta > 0$ are user-specified
regularization terms. 
Different from
existing works in~\cite{camisa2021distributed,falsone2017dual,vujanic2016decomposition}, 
we do not need to assume the uniqueness
of solutions.
Instead, $L_{\kappa}$ has a unique saddle point
due to being regularized, and this unique solution
will be found to approximately solve 
Problem~\ref{prob:matrix_Feasible_LP} without any additional
assumptions. 
Formally, we will approximately solve Problem \ref{prob:matrix_Feasible_LP}
by finding 
\begin{equation} \label{eq:saddledef}
(\hat{z}_{\kappa}, \hat{\lambda}_{\kappa}) = 
\argmin_{z \in Z}
\argmax_{\lambda \in \mathbb{R}^p_{+}}
\,\, L_{\kappa}(z, \lambda),
\end{equation}
which satisfies
\begin{equation} \label{eq:saddleineqs}
L_{\kappa}(\hat{z}_{\kappa}, \lambda) \leq 
L_{\kappa}(\hat{z}_{\kappa}, \hat{\lambda}_{\kappa}) \leq 
L_{\kappa}(z, \hat{\lambda}_{\kappa})
\end{equation}
for all~$z \in Z$ and~$\lambda \in \mathbb{R}^p_{+}$. 
However, due to the regularization, 
a saddle point of~$L_{\alpha,\delta}$
does not necessarily satisfy
the constraints in Problem~\ref{prob:matrix_Feasible_LP}. 

To guarantee constraint satisfaction, we will
bound the possible constraint violation of a saddle
point of~$L_{\alpha,\delta}$, and then we will tighten the constraints
by that amount. 

    
    \begin{lemma}[Constraint Tightening] \label{thm:newProblem}
        Let Assumptions 1 and 2 hold.
        Fix regularization parameters~$\alpha > 0$ and~$\delta > 0$. Then the maximum violation of any constraint
        in~$Az - b - \nu \leq 0$
        by the saddle point of~$L_{\alpha,\delta}$ 
        is denoted~$\phi$ and is bounded via
        \begin{equation} \label{eq:phivalue}
            \phi \leq \underset{1\leq j \leq m}{\max}\vert\vert A_{j} \vert \vert_{2} \cdot \left(\frac{c^{T}\bar{z}+\frac{\alpha}{2}\vert\vert\bar{z}\vert\vert^{2} + \vert\vert c\vert\vert\cdot r}{{\min\limits_{1\leq j\leq m}}-A_{j}\bar{z} + b_{j} + \nu_{j}}\right) \sqrt{\frac{\delta}{2\alpha}},
        \end{equation}
        where~$A_j$ is the~$j^{th}$ row of~$A$. 
    \end{lemma}
        \emph{Proof: } See Appendix~\ref{app:tightening}. \hfill $\blacksquare$
    
Now we can set the term $\phi$ 
equal to its upper bound and use it 
to enforce constraint satisfaction
by saddle points, even when~$L$ is regularized.
In Problem~\ref{prob:matrix_Feasible_LP}
we replace the constraint~$g_{\xi}(x, y) \leq 0$
with~$g_{\xi}(x, y) - \phi\bone \leq 0$,
which we assume to be feasible. 
Then we have
\begin{equation} \label{eq:RegularizedLagran}
 L_{\kappa}(z,\lambda)=c^{T}z+\frac{\alpha}{2}\vert\vert z\vert\vert^{2}+\lambda^{T}(Az-b-\nu-\phi\bone)-\frac{\delta}{2}\vert\vert\lambda\vert\vert^{2},
\end{equation}
where we use $\kappa=(\alpha,\delta,\phi)$ to simplify notation. 
One can find
a saddle point of the regularized Lagrangian
and then round the appropriate decision variables, 
and doing so gives
a point that satisfies all constraints 
of the original MILP. 

\begin{theorem} \label{thm:saddle_feasible}
Let Assumptions~\ref{as:Xcompact} and~\ref{as:slater} hold, and
consider the Lagrangian~$L_{\kappa}$
from~\eqref{eq:RegularizedLagran}. Let the unique saddle point
of~$L_{\kappa}$ be
denoted~$(\hat{z}_{\kappa}, \hat{\lambda}_{\kappa})$,
where~$\hat{z}_{\kappa} = (\hat{x}_{\kappa}, \hat{y}_{\kappa})$. 
Then the necessary rounding of its primal component, 
namely~$\round{\hat{z}_{\kappa}} = \big(\hat{x}_{\kappa}, \round{\hat{y}_{\kappa}}\big)$, satisfies
all constraints of Problem~\ref{prob:Original-MILP}. 
\end{theorem}
\emph{Proof: } 
By Lemma~\ref{thm:newProblem} we see that~$g_{\xi}\big(\hat{x}_{\kappa},\hat{y}_{\kappa}\big) \leq 0$, and thus
$\hat{z}_{\kappa} \in M_{\xi}$. 
By Proposition~\ref{prop:GranRounding}, rounding
any point in~$M_{\xi}$ produces a point that 
lies in~$M_{MILP}$. 
\hfill $\blacksquare$

\begin{remark} \label{rem:transcription}
Theorem~\ref{thm:saddle_feasible} shows that
we have faithfully transcribed Problem~\ref{prob:Original-MILP}
to a regularized saddle point problem, and it is this
regularized saddle point problem that is the focus of the
rest of the paper (it is formalized in Problem~\ref{prob:lksaddle}
below). 
Theorem~\ref{thm:saddle_feasible}
ensures that all constraints in Problem~\ref{prob:Original-MILP}
are satisfied by rounding a saddle
point, but it leaves open the extent to which such
a point is suboptimal. We next elaborate on the properties
of the regularized Lagrangian and then bound
this suboptimality. 
\end{remark}


As written, the regularized Lagrangian~$L_{\kappa}$
takes values in the domain~$Z \times \mathbb{R}^p_{+}$, which is 
unbounded. To facilitate the forthcoming analysis and algorithm
design, we compute a compact set that contains~$\hat{\lambda}_{\kappa}$.
Since the set~$Z$ is compact, it is in particular bounded,
and therefore there exists 
\begin{equation} \label{eq:r}
r = \max_{z \in Z} \|z\|_2,
\end{equation}
and a ball of radius~$r$ at the origin that
contains~$Z$, namely
\begin{equation} \label{eq:b0}
Z \subseteq B_0(r) = \{z \in \mathbb{R}^{n+m} : \|z\|_2 \leq r\},
\end{equation}
which we use in the next result. 
\begin{lemma} \label{lem:saddlePoint}
Let Assumptions~\ref{as:Xcompact} and~\ref{as:slater} hold. Then
\begin{equation}
\hat{\lambda}_{\kappa} \in \Lambda = 
\left\{ \lambda\in\mathbb{R}_{+}^{p}:
\|\lambda\|_1 \leq
\frac{c^{T}\bar{z}+\frac{\alpha}{2}\vert\vert\bar{z}\vert\vert^{2} + \vert\vert c\vert\vert\cdot r}{\min\limits_{1\leq j\leq m} 
-A_{j}\bar{z} + b_{j} + \nu_{j} + \phi}\right\} ,
\end{equation}
where~$c$ is from~\eqref{eq:c}, $\bar{z}$ is from~\eqref{eq:zbar},
$A$ is from~\eqref{eq:A} and~$A_j$ denotes the~$j^{th}$
row of~$A$ for~$j \in \{1, \ldots, p\}$,
$b$ is from~\eqref{eq:b},
$\nu$ is from~\eqref{eq:nu},
$r$ is from~\eqref{eq:r}, $\phi$ is from~\eqref{eq:phivalue},
and~$\alpha$ is the primal regularization term in~$L_{\kappa}$
from~\eqref{eq:RegularizedLagran}. 
\end{lemma}
\emph{Proof: }
See Appendix~\ref{ss:lem1}. \hfill $\blacksquare$

Now that we have a compact set that contains $\hat{\lambda}_{k},$ we formally state the saddle point problem that we solve.

\begin{problem} \label{prob:lksaddle}
Find the unique saddle point~$(\hat{z}_{\kappa}, \hat{\lambda}_{\kappa})$, where
\begin{equation}
(\hat{z}_{\kappa},\hat{\lambda}_{\kappa})=
\argmin_{z \in Z}
\argmax_{\lambda \in \Lambda} L_{\kappa}(z, \lambda)
\end{equation}
and where~$L_{\kappa}$ is from~\eqref{eq:RegularizedLagran}. 
\end{problem}

\subsection{Suboptimality Bounds}
We now bound the suboptimality induced by the regularization of $L$. 

\begin{lemma}
\label{lem:regError} Let Assumptions 1 and 2 hold.
Then
\begin{equation}
    \vert c^{T}\hat{z}_{\kappa}-c^{T}z^{*}\vert \leq \vert\vert c\vert\vert \left(\frac{c^{T}\bar{z}+\frac{\alpha}{2}\vert\vert\bar{z}\vert\vert^{2} + \vert\vert c\vert\vert r}{\min\limits_{1\leq j\leq m} -A_{j}\bar{z} + b_{j} + \nu_{j} + \phi}\right)\sqrt{\frac{\delta}{2\alpha}}+\frac{\alpha}{2}r,
\end{equation}
where $\hat{z}_{\kappa}$ is the solution to 
Problem~\ref{prob:lksaddle}, $z^{*}$ is the 
solution to Problem~\ref{prob:matrix_Feasible_LP},
$\bar{z}$ is from Assumption~\ref{as:slater}, $\alpha$ is the primal regularization
term, $\delta$ is the dual regularization term, $r$ is from
 (\ref{eq:r}), and $\phi$ is from~\eqref{eq:phivalue}. 
\end{lemma}
\emph{Proof: } See Appendix~\ref{app:regError}. \hfill $\blacksquare$

We will next work to bound 
the difference in optimal costs between Problem
\ref{prob:Original-MILP} and Problem~\ref{prob:lksaddle}. 
Toward doing so, 
we next derive suboptimality bounds 
between the solution to Problem~\ref{prob:matrix_Feasible_LP} and the solution to Problem \ref{prob:Original-MILP}. 
These bounds are in terms of a Hoffman constant~\cite{aze2002sensitivity}. 

\begin{definition}[Hoffman constants] \label{def:hoffman}
Let~$C \in \mathbb{R}^{p \times n}$, $D \in \mathbb{R}^{p \times m}$,
and~$b \in \mathbb{R}^p$, and suppose that
the polyhedron
\begin{equation}
P_{C,D,b} = \{(x, y) \in \mathbb{R}^n \times \mathbb{R}^m \mid Cx + Dy \leq b\}
\end{equation}
is non-empty. Then there exists a constant~$\sigma(C,D)$ such that,
for all~$(u, v) \in \mathbb{R}^n \times \mathbb{R}^m$ we have
\begin{equation} \label{eq:hoffdef}
\min_{(x,y) \in P_{C,D,b}} \left\|\left(\begin{array}{c} x \\ y \end{array}\right)
- \left(\begin{array}{c} u \\ v \end{array}\right)\right\|_2 \leq \sigma(C,D) \|(Cu + Dv - b)_{+}\|_2,
\end{equation}
where
\begin{equation}
(Cu + Dv -b)_{+} = \max\{0, Cu + Dv - b\}, 
\end{equation}
with the maximum applied
componentwise. Any~$\sigma$ that satisfies~\eqref{eq:hoffdef} is
\emph{a Hoffman constant} for the constraints~$Cx + Dy \leq b$. 
\end{definition}

The following result bounds
the difference in costs
between the solutions to Problem~\ref{prob:Original-MILP} and Problem~\ref{prob:matrix_Feasible_LP}.

\begin{lemma}\label{lem:HoffmanError}
Let $\sigma$
denote a Hoffman constant of the system of inequalities $Ex + Fy \leq h$,
let~$z^*_{MILP}$ denote the solution to Problem~\ref{prob:Original-MILP},
and let~$z^* = (x^{*,T}, y^{*,T})^T$ denote
the solution to Problem~\ref{prob:matrix_Feasible_LP}.
Then 
\begin{equation}
    \vert c^{T}z_{\text{MILP}}^{*}-c^{T}z^*\vert \leq
    \frac{\|c\|_2}{2}\big(\sigma \|F\|_{\infty} + 1 \big).
\end{equation}
\end{lemma}
\emph{Proof:} See Proposition~4.2 in~\cite{stein2016error}. \hfill $\blacksquare$

Next, the following theorem relates
the optimal cost in Problem~\ref{prob:Original-MILP},
which is the original MILP problem, to the optimal
cost in Problem~\ref{prob:lksaddle}, which is the regularized 
saddle point problem that we solve. 
This is the main error bound that we need, and it accounts
for the effects of both relaxing Problem~\ref{prob:Original-MILP}
and regularizing its associated Lagrangian. 

\begin{theorem}[Suboptimality Error Bound] \label{thm:suboptimalityEB}
    Let~$z^*_{MILP}$ be the solution to Problem~\ref{prob:Original-MILP}
    and let~$\hat{z}_{\kappa}$ be the solution to Problem~\ref{prob:lksaddle}.
    Let $\sigma$ be a Hoffman constant 
    for the constraints~$Ex + Fy \leq h$. 
    Then
    \begin{equation}
        \vert c^{T}z_{\text{MILP}}^{*}-c^{T}\round{\hat{z}_{\kappa}}\vert \leq  \frac{\vert\vert c\vert\vert_{2}}{2}\big(\sigma\vert\vert F\vert\vert_{\infty}+1\big) + \frac{\vert\vert d\vert\vert_{1}}{2} + \vert\vert c\vert\vert\left(\frac{c^{T}\bar{z}+\frac{\alpha}{2}\vert\vert\bar{z}\vert\vert^{2} + \vert\vert c\vert\vert\cdot r}{\underset{1\leq j\leq m}{\min}-A_{j}\bar{z}+b_{j}+\nu_{j}+\phi}\right)\sqrt{\frac{\delta}{2\alpha}}+\frac{\alpha}{2} r,
    \end{equation}
    where~$\round{\hat{z}_{\kappa}} = (\hat{x}_{\kappa}^T, \round{\hat{y}_{\kappa}}^T)^T$. 
    \end{theorem}
    \emph{Proof: } 
     Applying the triangle inequality gives
    \begin{equation}
    \big\vert(c^{T}z_{\text{MILP}}^{*}-c^{T}z^*)+(c^{T}z^{*}-c^{T}\hat{z}_{\kappa}) + 
    (c^T\hat{z}_{\kappa} - c^T\round{\hat{z}_{\kappa}})
    \big\vert 
    \leq\vert c^{T}z_{\text{MILP}}^{*}-c^{T}z^*\vert + \vert c^{T}z^{*}-c^{T}\hat{z}_{\kappa}\vert + |c^T\hat{z}_{\kappa} - c^T\round{\hat{z}_{\kappa}}|.
    \end{equation}
    We conclude by applying
    Lemmas~\ref{lem:regError} and~\ref{lem:HoffmanError}, and  
    noting that
    \begin{equation}
    |c^T\round{\hat{z}_{\kappa}} - c^T\hat{z}_{\kappa}| = |d^T(\round{\hat{y}_{\kappa}} - \hat{y}_{\kappa})| \leq \|d\|_1\|\round{\hat{y}_{\kappa}} - \hat{y}_{\kappa}\|_{\infty} \leq \frac{\|d\|_1}{2}
    \end{equation}
    which follows from H{\"o}lder's Inequality and the fact that no entry of~$\hat{y}_{\kappa}$ changes by more than~$1/2$
    when it is rounded. 
    \hfill $\blacksquare$
   
    Theorem~\ref{thm:suboptimalityEB} shows that the 
    various modifications made to Problem~\ref{prob:Original-MILP}
    result in bounded, quantifiable sub-optimality, and this bound
    can be applied by users to select the regularization
    parameters~$\alpha$ and~$\delta$ in order to control
    this error. Theorem~\ref{thm:suboptimalityEB}
    is a worst-case bound, and we will show in
    Section~\ref{sec:Simulation} that realistic problems
    can result in negligible sub-optimality.

\section{Distributed Algorithm\label{sec:Distributed-Algorithm}}
In this section we define the main algorithm that we use to
find the saddle point of~$L_{\kappa}$, 
analyze its convergence, and describe how its final result can
be mapped back to an approximate solution to Problem~\ref{prob:Original-MILP}. 

\subsection{Algorithm Definition}
The distributed algorithm that we develop consists of three types of operations: 
(i) primal variable computations, 
(ii) dual variable computations, and (iii)
communications of the results of these computations. 
Recall from Problem~\ref{prob:lksaddle} that~$z=(x^T, y^T)^T \in \mathbb{R}^{m+n}$. 
Suppose that there are~$N = N_1 + N_2$ agents
to find the saddle point~$(\hat{z}_{\kappa}, \hat{\lambda}_{\kappa})$ of~$L_{\kappa}$, 
of which there are (i)~$N_1$ primal agents 
(agents that perform computations on primal decision
variables)
indexed over the set $A_{P}=[N_1]$, 
and (ii) $N_2$ dual agents (agents that perform computations on dual decision
variables), indexed
over the set $A_{D}=[N_2]$.
Each primal and dual agent has partial knowledge of
Problem \ref{prob:matrix_Feasible_LP} and only updates its assigned block
of the decision variables.

We partition $z$ into $N_1$ blocks due to the separability of the objective function in Problem~\ref{prob:matrix_Feasible_LP}. Formally, 
\begin{equation}
z = (z_{[1]}^T, \ldots, z_{[N_1]}^T)^T \in \mathbb{R}^{m+n},
\end{equation}
where~$z_{[i]} \in \mathbb{R}^{n_i}$ for~$i \in A_P$ and~$\sum_{i \in A_P} n_i = m+n$. 
We emphasize that~$z_{[i]}$ does not need to be a scalar. By combining Assumption~\ref{as:Xcompact}
and Remark~\ref{rem:yprod}, we see that
\begin{equation}
Z = X_1 \times X_2 \times \cdots \times X_n \times Y_1 \times \cdots \times Y_m,
\end{equation}
and the constraint set that contains~$z_{[i]}$ is the Cartesian product of the~$X_i$'s and~$Y_i$'s
that contain its entries. We denote this set by~$Z_i$, and thus
the set constraint 
\begin{equation}
z_{[i]} \in Z_i \subseteq \mathbb{R}^{n_i}
\end{equation}
is enforced by agent~$i$. 

We partition the dual 
variable~$\lambda \in \mathbb{R}^p_{+}$ into~$N_2$ blocks via
\begin{equation}
\lambda = (\lambda_{[1]}^T, \ldots, \lambda_{[N_2]}^T)^T,
\end{equation}
where~$\lambda_{[q]} \in \mathbb{R}^{r_q}_{+}$
for~$q \in A_D$ and~$\sum_{q \in A_D} r_q = p$. 
Each block~$\lambda_{[q]}$ is constrained to lie in the set
\begin{equation} \label{eq:lambdaqbound}
\Lambda_q \!=\! \left\{\!\lambda_{[q]} \in \mathbb{R}_{+}^{r_q} \!:\! \vert\vert\lambda_{[q]}\vert\vert_{1} \!\leq\! \frac{c^{T}\bar{z}+\frac{\alpha}{2}\vert\vert\bar{z}\vert\vert^{2} + \vert\vert c\vert\vert r}{\underset{1\leq j\leq m}{\min} \!\!{-A}_{j}\bar{z}+b_{j}+\nu_{j}+\phi}\right\}.
\end{equation}

The upper bound in the definition of~$\Lambda_q$ is the same as that in~$\Lambda$ in Lemma~\ref{lem:saddlePoint} and
thus is a bound on the~$1$-norm of 
each block of~$\lambda$ as well. 
The bound in~\eqref{eq:lambdaqbound} enlarges the dual feasible region beyond~$\Lambda$, though
it still gives a compact feasible region that contains the dual optimum, which is what
we require to solve Problem~\ref{prob:lksaddle}.

Computations and communications are asynchronous, which leads to disagreements
among the values of the decision variables that agents store onboard.
Therefore, there is a need to 
compute the saddle point~$(\hat{z}_{\kappa}, \hat{\lambda}_{\kappa})$
with an update law that is robust to asynchrony. We use a decentralized
form of the classic Uzawa iteration, also called gradient descent-ascent \cite{arrow1958studies},
to find a saddle point of the regularized Lagrangian~$L_{\kappa}$ associated with
Problem \ref{prob:matrix_Feasible_LP}. 

Each primal agent updates the block of the primal variables it stores onboard. To do this, 
primal agent $i$ stores onboard
itself a local copy of its primal variables, 
denoted $z^{i}_{[i]} \in Z_i$,
and a local copy of the vector of dual variables, 
denoted $\lambda^{i} \in \Lambda$. 
Similarly, dual agent $q$ stores a local
primal vector onboard, denoted $z^{q} \in Z$,
and its block of the dual vector, denoted $\lambda^{q}_{[q]} \in \Lambda_q$. 
Each dual agent updates
its block of the dual variables in the copy that it stores onboard.
Each dual variable corresponds to an inequality constraint on the primal
variables, and each dual agent sends values
of its updated block to all primal
agents whose decision variables are in the constraints that correspond
to those dual variables.

We define $D\subseteq\mathbb{N}$ as the set of times at which all dual agents
compute updates to their decision variables, 
and we define~$K^{i}\subseteq\mathbb{N}$ as the set of times
at which primal agent $i\in A_{P}$ computes updates
to its decision variables. To state the algorithm, we use
$k$ as the iteration counter used by the primal agents
and $t$ as the iteration counter used
by dual agents. 
We emphasize that the sets $K^{i}$ and $D$ are tools for analysis and discussion, and they need not be known by 
any agent. 
Additionally,
$\lambda^{i}(k)$ denotes the value of~$\lambda$ stored onboard primal
agent~$i$ at time~$k$, and~$z^q(k)$ denotes the value of~$z$
stored onboard dual agent~$q$ at time~$k$. 

We make the following assumption about delays, 
termed ``partial asynchrony''~\cite{bertsekas2015parallel}. 

\begin{assumption} \label{as:partialasynchrony}
Let $K^{i}$ be the set of times at which primal agent~$i$ performs
a computation. 
Then there exists~$B \in \mathbb{N}$ such that 
for every $i\in A_{P}$, at least one of the elements of the
set $\{k,k+1,...,k+B-1\}$ is in $K^{i}$. 
\end{assumption}
This assumption ensures that
all primal agents perform at least one computation every~$B$ timesteps.

Primal agent $i$ updates its block~$z^i_{[i]}$ via projected gradient descent
at each time $k\in K^{i}$. If $k\notin K^{i}$,
then agent $i$ does not update and~$z^i_{[i]}$ is held constant. 
At each timestep $t\in D$, dual
agent~$q$ updates via projected gradient ascent onto the set $\Lambda_{q}$
from~\eqref{eq:lambdaqbound}.

Primal agents $i$'s computations take the form
\begin{equation}
z_{[i]}^{i}(k+1) \!=\! \begin{cases}
\Pi_{Z_{i}}[z_{[i]}^{i}(k) \!-\! \gamma\nabla_{z_{[i]}}L_{\kappa}(z^{i}(k),\lambda^{i}(k))] & k \!\in\! K^{i}\\
z_{i}^{i}(k) & k \!\not\in\! K^i
\end{cases}
\end{equation}
where~$\nabla_{z_{[i]}} = \frac{\partial}{\partial z_{[i]}}$. 

For all~$q \in A_D$, dual agent~$q$ only performs an update
occasionally. Specifically, it only performs an update after
the counter~$k$ has increased by some amount that is divisible
by~$B$. 
When dual agent~$q$ updates, 
it computes
\begin{equation}
\lambda_{[q]}^{q}(tB) = 
\Pi_{\Lambda_{q}}[\lambda_{[q]}^{q}(tB-1)+\beta\nabla_{\lambda_{[q]}}L_{\kappa}(z^{q}(tB),\lambda^{q}(tB-1))],
\end{equation}
and between updates it holds the value of~$\lambda^q_{[q]}$ constant. 
After updating, dual agent~$q$ sends the new value of~$\lambda^{[q]}_q$ to
all primal agents that need it in their computations. 

\subsection{Communications and Algorithm Statement}
A primal agent only needs to receive 
communcations of updated values of 
a dual block if
that primal agent's decision variables appear in the constraints
corresponding to those dual variables. Similarly, a dual agent
only needs to receive communications from a primal agent if
that primal agent's decision variables appear in the constraints
encoded in the dual variables that the dual agent updates. 
In contrast to related work in the literature~\cite{ubl19,hochhaus18},
there is no communication between pairs of primal agents, and 
similarly there is no communication
between pairs of dual agents. 
We formalize these notions in the following definition.

\begin{definition}[Essential Neighbors] \label{def:essneighbors}
Let~$i \in A_P$ be the index of a primal agent, and let~$q \in A_D$ be the index of a dual
agent. Suppose that agent~$i$ updates the entries of~$z$ in the block~$z^i_{[i]}$,
and suppose that agent~$q$ updates the entries of~$\lambda$ in the block~$\lambda^q_{[q]}$.
Then agents~$i$ and~$q$ communicate if and only if at least one entry of~$z_{[i]}$
appears in at least one of the constraints~$g_{[q]}$. In such a case, we say that
agent~$i$ is an ``essential neighbor'' of agent~$q$ and that agent~$q$
is an ``essential neighbor'' of agent~$i$. 
\end{definition}

The following example illustrates this idea. 

\begin{example} \label{ex:comms}
Consider the problem
\begin{align}
&\underset{x \in \mathbb{R}^2, y \in \mathbb{R}^2}{\textnormal{minimize }} \,\, x_1 + x_2 + y_1 + y_2 \\
&g(x,y) = \left(\begin{array}{c} g_1(x, y) \\ g_2(x, y)\end{array}\right) = \left(\begin{array}{c}
x_1 + y_1 + 1 \\
x_2 + y_2 + 1
\end{array}\right) \leq 0.
\end{align}
Suppose there are two primal agents and two dual agents used to solved this problem.

Suppose that primal agent~$1$ updates~$x^1_{[1]}$ and~$y^1_{[1]}$ and
that primal agent~$2$ updates~$x^2_{[2]}$ and~$y^2_{[2]}$. 
Suppose that dual agent~$1$ updates~$\lambda^1_{[1]}$, which
corresponds to~$g_1$, and suppose that dual agent~$2$ updates~$\lambda^2_{[2]}$, which corresponds to~$g_2$. 

\begin{figure}
\centering
\includegraphics[width=3.3in,,draft=false]{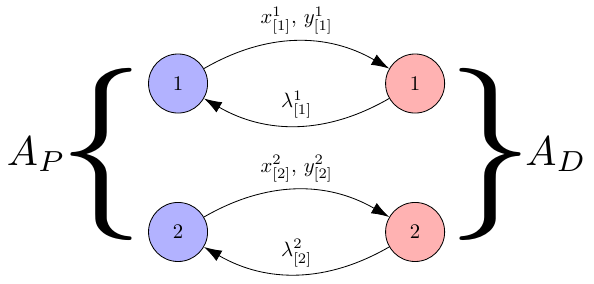}
\caption{Agents' essential neighbors in Example~\ref{ex:comms}.}
\label{fig:ex1}
\end{figure}

Then communications with essential neighbors as defined in Definition~\ref{def:essneighbors} produces
the communication topology shown in Figure~\ref{fig:ex1}. There are no
communications between pairs of primal agents and there are no communications between pairs of dual agents. The communications
between primal and dual agents are only between essential neighbors, and, in particular,
they are not all-to-all. 
\end{example}

It is required that all primal agents use the same values
of the dual blocks in their computations. This has been shown
to be necessary
for convergence~\cite{hendrickson2021totally}; any mechanism to
enforce this agreement can be used. 

We present the full parallelized block-based saddle point algorithm
in Algorithm \ref{alg:algorithm1}. 

\begin{algorithm}[ht] 
\caption{Distributed saddle point finding algorithm}
\label{alg:algorithm1}
\SetAlgoLined
\KwIn{
All data from Problem~\ref{prob:lksaddle}, 
primal vectors $z^i_{[i]}(0)\in Z$ for all~$i \in A_P$, 
dual vectors~$\lambda^q_{[q]}(0) \in \Lambda$ for all~$q \in A_D$, 
stepsizes~$\gamma > 0$ and~$\beta > 0$, regularization parameters~$\alpha > 0$
and~$\delta > 0$} 
\KwOut{The saddle point~$\hat{z}_{\kappa}$}
\For{$k = 0, 1, 2,$}{
\For{$i \in A_P$}{
\If{$k\in K^{i}$}{
{$z_{[i]}^{i}(k+1)=\Pi_{Z_{i}}[z_{[i]}^{i}(k)-\gamma\nabla_{z_{[i]}}L_{\kappa}(z^{i}(k),\lambda^{i}(k))]$} 
}
\Else{
$z_{[i]}^{i}(k+1)=z_{[i]}^{i}(k)$
}
}
\For{$q \in A_D$}{
\If{$k=tB\text{ with }t\in D$}{
\For{$i \in A_P$}{Agent~$i$ sends~$z^i_{[i]}$ to its essential neighbors}
}
\For{$q\in A_{D}$}{
$\lambda_{[q]}^{q}(tB)=\Pi_{\Lambda_{q}}[\lambda_{[q]}^{q}(tB-1)+\beta\nabla_{\lambda_{[q]}}L_{\kappa}(z^{q}(tB),\lambda^{q}(tB-1))]$ \\
Agent~$q$ sends~$\lambda^{q}_{[q]}$ to its essential neighbors
}

}
}
\For{$i \in A_P$}{
If agent~$i$ was responsible for computing any entries
of~$y$, then it rounds those entries to the nearest integer (ties are broken arbitrarily)
}
\end{algorithm}

\subsection{Convergence Analysis}
We next need to verify that Algorithm \ref{alg:algorithm1} will converge
to the saddle point of $L_{\kappa}$. We start with a proof of convergence
for the primal agents when they have a fixed dual variable onboard.
In it, we define the primal agents' true iterate for all~$k \in \mathbb{N}$ as
\begin{equation}
z(k) = \big(z^1_{[1]}(k)^T, z^2_{[2]}(k)^T, \ldots, z^{N_1}_{[N_1]}(k)^T\big). 
\end{equation}

\begin{lemma}[Primal Convergence Between Dual Updates] \label{lem:primalConvergence} Let Assumptions
\ref{as:Xcompact}-\ref{as:partialasynchrony} hold, and consider using Algorithm \ref{alg:algorithm1} 
to solve Problem~\ref{prob:lksaddle}. Fix $t\in D$ and let $t'$ be the
smallest element of D that is greater than t. Fix $\lambda(tB)\in \Lambda$
and define $\hat{z}(tB)=\argmin_{z\in Z}L_{\kappa}(z,\lambda(tB))$
as the point that the primal agents would converge to with $\lambda(tB)$
held fixed. For agents executing Algorithm \ref{alg:algorithm1} at
times $k\in\{tB,tB+1,\ldots,t'B\},$ there exists a scalar $\gamma_{1}>0$
such that for~$\gamma \in (0, \gamma_1)$ we have
\begin{equation}
\vert\vert z(t'B)-\hat{z}(tB)\vert\vert\leq(1-\theta\gamma)^{t'-t}\vert\vert z(tB)-\hat{z}(tB)\vert\vert,
\end{equation}
where $\theta$ is a positive constant, and $(1-\theta\gamma) \in [0,1)$.
\end{lemma}
\emph{Proof: }
We show that Assumptions A and B in \cite{tseng1991rate} are satisfied,
which enables the use of Proposition 2.2 from the same reference to show
convergence.
Assumption A requires that 
\begin{enumerate} 
\item $L_{\kappa}(\cdot,\lambda(tB))$ is bounded from below on $Z$, 
\item the set $Z$ contains at least one point $v$ such that $v=\Pi_{Z}[v-\nabla_{z}L_{\kappa}(v,\lambda(tB))]$,
\item $\nabla_{z}L_{\kappa}(\cdot,\lambda(tB))$ is Lipschitz on $Z$. 
\end{enumerate}

First, we see that for all~$z \in Z$ we have 
\begin{equation}
\nabla_z L_{\kappa}\big(z, \lambda(tB)\big) = c + \alpha z
+ A^T\lambda(tB),
\end{equation}
which is Lipschitz in~$z$ with constant~$\alpha$. Then Assumption A.3 is satisfied. 


Next, for every choice of $\lambda(tB)\in \Lambda$, the function $L_{\kappa}(\cdot,\lambda(tB))$
is bounded from below because it is continuous and its domain~$Z$ is compact. 
Then Assumption A.2 is satisfied. For every $t\in\mathbb{N}$ and
every fixed $\lambda(tB)\in \Lambda$, the strong convexity of $L_{\kappa}(\cdot,\lambda(tB))$
implies that is has a unique minimum over $Z$, denoted $\hat{z}(tB)$.
This point is 
the unique fixed point of the projected gradient descent mapping
and Assumption A.1 is satisfied. Then all conditions of Assumption
A in \cite{tseng1991rate} are satisfied. 

Assumption B in \cite{tseng1991rate} requires 
\begin{enumerate}
\item the isocost curves of $L_{\kappa}(\cdot,\lambda(tB))$ to be separated,
and 
\item the ``error bound condition'', stated as (2.5) in \cite{tseng1991rate},
is satisfied.
\end{enumerate}

It is observed in~\cite{tseng1991rate} 
that both criteria are satisfied by functions that are strongly convex
with Lipschitz gradients. In this work,  $L_{\kappa}(\cdot,\lambda(tB))$ has
both of these properties. Then Assumption B is satisfied as well, and 
an application of Proposition~2.2 in \cite{tseng1991rate} completes the proof. \hfill $\blacksquare$

Next, we derive a proof of convergence for dual agents across two consecutive timesteps.

\begin{lemma}[Dual Convergence Between Two Time Steps] \label{lem:(Dual-Convergence-Between} 
Let Assumptions \ref{as:Xcompact}-\ref{as:partialasynchrony} 
hold and consider the use of Algorithm
\ref{alg:algorithm1} to solve Problem~\ref{prob:lksaddle}. Let
the dual step size be $0<\beta<\min\left\{\frac{2\alpha}{\vert\vert A \vert\vert_{2}+2\alpha\delta},\frac{2\delta}{1+\delta^{2}}\right\}$, 
let the primal stepsize be~$\gamma \in (0, \gamma_1)$, 
and let $t_{1}$ and $t_{2}$ denote two consecutive times at which dual
updates have occurred, with $t_{1}<t_{2}$ and $t_{1},t_{2}\in D.$
Then Algorithm~1 produces dual variables~$\lambda(t_1B)$
and~$\lambda(t_2B)$ that satisfy 
\begin{equation}
\vert\vert\lambda(t_{2}B)-\hat{\lambda}_{\kappa}\vert\vert^{2}\leq q_{d}\vert\vert\lambda(t_{1}B)-\hat{\lambda}_{\kappa}\vert\vert^{2} + 4r^2q_{d}q_{p}^{2(t_{2}-t_{1})}\vert\vert A\vert\vert_{2}^{2} +8r^2\beta^{2}q_{p}^{t_{2}-t_{1}}\vert\vert A\vert\vert_{2}^{2},
\end{equation}
where $q_{d} = (1-\beta\delta)^{2}+\beta^{2}\in[0,1)$, $q_{p} = (1-\theta\gamma)\in[0,1)$, 
$r$ is from~\eqref{eq:r}, and $\hat{\lambda}_{\kappa}$ is
the dual component of the unique saddle point of $L_{\kappa}$. 
\end{lemma}
\emph{Proof: }
See Appendix~\ref{ss:lem7}.
\hfill $\blacksquare$

In Lemma~\ref{lem:(Dual-Convergence-Between}, the terms that do not explicitly depend on the dual variable
are dependent on regularization parameters and the number of primal agent iterations
that occur between dual agents' updates. Now we will use the convergence
rate for two consecutive timesteps and derive a proof of convergence
for the dual variables to the optimum.

\begin{theorem}[Dual Convergence to Optimum] \label{thm:dualConvergence}
Let
Assumptions \ref{as:Xcompact}-\ref{as:partialasynchrony} hold and consider the use of Algorithm \ref{alg:algorithm1}
to solve Problem~\ref{prob:lksaddle}. Let the dual step size
be $0<\beta<\min\left\{\frac{2\alpha}{\vert\vert A\vert\vert_{2}+2\alpha\delta},\frac{2\delta}{1+\delta^{2}}\right\}$
and let the primal stepsize be~$\gamma \in (0, \gamma_1)$. 
Additionally, let $t_{n}$ denote the $n^{\text{th}}$ entry in D
where $t_{1}<t_{2}< \cdots <t_{n}$, which means that $t_{n}B$ is the time
when the $n^{\text{th}}$ dual update occurs for all dual agents.
Then the dual agents executing Algorithm \ref{alg:algorithm1} generate
dual variables that satisfy
\begin{equation}
 \vert\vert\lambda\left(t_{n}B\right)-\hat{\lambda}_{\kappa}\vert\vert^{2}  \leq q_{d}^{n}\vert\vert\lambda(0)-\hat{\lambda}_{\kappa}\vert\vert^{2}  
 +\left(4r^2q_{d}q_{p}^{2}\vert\vert A\vert\vert_{2}^{2}+8r^2\beta^{2}q_{p}\vert\vert A\vert\vert_{2}^{2}\right)\sum_{i=0}^{n}q_{d}^{i},
\end{equation}
where $q_{d}=(1-\beta\delta)^{2}+\beta^{2} \in [0, 1),$ $q_{p}=(1-\theta\gamma)\in[0,1)$, $r$ is from~\eqref{eq:r}, 
and $\hat{\lambda}_{\kappa}$ is
the dual component of the unique saddle point of $L_{\kappa}$. 
\end{theorem}
\emph{Proof: }
We apply Lemma \ref{lem:(Dual-Convergence-Between} twice to find 
\begin{align}
\vert\vert\lambda\left(t_{n}B\right)-\hat{\lambda}_{\kappa}\vert\vert^{2} & \leq q_{d}\vert\vert\lambda\left(t_{n-1}B\right)-\hat{\lambda}_{\kappa}\vert\vert^{2} +4r^2q_{d}q_{p}^{2\left(t_{n}-t_{n-1}\right)}\vert\vert A\vert\vert_{2}^{2} +8r^2\beta^{2}q_{p}^{t_{n}-t_{n-1}}\vert\vert A\vert\vert_{2}^{2}\\
 & \leq q_{d}^{2}\vert\vert\lambda\left(t_{n-2}B\right)-\hat{\lambda}_{\kappa}\vert\vert^{2} +4r^2q_{d}^{2}q_{p}^{2}\vert\vert A\vert\vert_{2}^{2}+8r^2\beta^{2}q_{d}q_{p}\vert\vert A\vert\vert_{2}^{2} + 4r^2q_{d}q_{p}^{2}\vert\vert A\vert\vert_{2}^{2}+8r^2\beta^{2}q_{p}\vert\vert A\vert\vert_{2}^{2},
 \end{align}
 where the second inequality follows because $t_{n}-t_{n-1}\geq1$.
 Recursively applying Lemma~\ref{lem:(Dual-Convergence-Between} then gives
 the result. 
\hfill $\blacksquare$

Finally this bound is used to formulate a convergence rate for
the primal variables in Algorithm~\ref{alg:algorithm1}. 

\begin{theorem}[Overall Primal Convergence] \label{thm:overallConvergence} 
Let Assumptions \ref{as:Xcompact}-\ref{as:partialasynchrony}
hold, let the dual stepsize satisfy $0<\beta<\min\left\{\frac{2\alpha}{\vert\vert A \vert \vert_{2}+2\alpha\delta},\frac{2\delta}{1+\delta^{2}}\right\}$, and let the primal stepsize
satisfy~$\gamma \in (0, \gamma_1)$. 
Consider using Algorithm~\ref{alg:algorithm1} to solve Problem~\ref{prob:lksaddle}, 
and let $t_{n}$ denote
the $n^{\text{th}}$ entry in D where $t_{1}<t_{2}< \cdots <t_{n}$.
That is, $t_{n}B$ is the time at which the $n^{\text{th}}$ dual
update occurs across all dual agents. Then primal agents executing Algorithm
\ref{alg:algorithm1} generate primal iterates that satisfy
\begin{equation}
\vert\vert z(t_{n}B)-\hat{z}_{k}\vert\vert\leq 2q_{p}^{t_{n}-t_{n-1}}r +\frac{\vert\vert A\vert\vert_{2}}{\alpha}\left(q_{d}^{n-1}\vert\vert\lambda(0)-\hat{\lambda}_{\kappa}\vert\vert^{2} + \left(4r^2q_{d}q_{p}^{2}\vert\vert A\vert\vert_{2}^{2}+8r^2\beta^{2}q_{p}\vert\vert A\vert\vert_{2}^{2}\right)\sum_{i=0}^{n-1}q_{d}^{i}\right)^{1/2},
\end{equation}
where $q_{p}=(1-\theta\gamma)\in[0,1),$ $q_{d}$$=(1-\beta\delta)^{2}+\beta^{2}\in[0,1)$\emph{,
and $\hat{\lambda}_{\kappa}$ is the dual at the optimum of $L_{\kappa}.$}
\end{theorem}
\emph{Proof: } We use
$\hat{z}(t_iB) = \argmin_{z \in Z} L_{\kappa}\big(z, \lambda(t_iB)\big)$
for all~$i \in [n]$. 
Using the triangle inequality, we have
\begin{equation}
\|z(t_nB) - \hat{z}_{\kappa}\| = \|z(t_nB) - \hat{z}(t_{n-1}B) + \hat{z}(t_{n-1}B) - \hat{z}_{\kappa}\| 
\leq \|z(t_nB) - \hat{z}(t_{n-1}B)\|
+ \|\hat{z}(t_{n-1}B) - \hat{z}_{\kappa}\|. \label{eq:t3tri} 
\end{equation}

From Lemma~\ref{lem:primalConvergence}, we see that
\begin{equation}
\|z(t_nB) - \hat{z}(t_{n-1}B)\| \leq (1 - \theta\gamma)^{t_n-t_{n-1}}
\|z(t_{n-1}B) - \hat{z}(t_{n-1}B)\| 
    \leq 2q_p^{t_n - t_{n-1}}r, \label{eq:t3qp}
\end{equation}
where the second inequality follows from the definition of~$r$ in~\eqref{eq:r}
Next, by~\cite[Lemma 4.1]{koshal2011multiuser}, we see that 
\begin{equation}
\|\hat{z}(t_{n-1}B) - \hat{z}_{\kappa}\| \leq
\frac{\|A\|_2}{\alpha}\|\lambda(t_{n-1}B) - \hat{\lambda}_{\kappa}\|. \label{eq:t3dualbound}
\end{equation}
Using~\eqref{eq:t3qp} and~\eqref{eq:t3dualbound}
in~\eqref{eq:t3tri} gives
\begin{equation}
\|z(t_nB) - \hat{z}_{\kappa}\| \leq 2q_p^{t_n-t_{n-1}}r
+ \frac{\|A\|_2}{\alpha}\|\lambda(t_{n-1}B) - \hat{\lambda}_{\kappa}\|,
\end{equation}
and bounding the last term using Theorem~\ref{thm:dualConvergence} completes the proof. 
\hfill $\blacksquare$

\begin{remark}
The upper bound in Theorem~\ref{thm:overallConvergence} is non-zero in general, which
means that primal agents do not exactly 
compute~$\hat{z}_{\kappa}$ in finite time.
Therefore, while Theorem~\ref{thm:saddle_feasible} guarantees that
rounding~$\hat{z}_{\kappa}$ produces
a point that satisfies the constraints in Problem~\ref{prob:Original-MILP}, 
rounding agents' final
primal iterate may not. To account for this, one could simply increase the value of~$\phi$
in~\eqref{eq:phivalue} to tighten the constraints there. However, in practice the 
distance~$\|z(tB) - \hat{z}_{\kappa}\|$ is often negligible, which means that tightening
the constraints further is not necessary to attain a feasible point after rounding. 
We illustrate this behavior in simulations in the next section. 
\end{remark}
\section{Numerical Results\label{sec:Simulation}}
In this section, we solve the generalized assignment problem (GAP)~\cite[Chapter 2]{conforti2014integer} to illustrate
our developments in a practical application. 
Among other applications, GAP can be applied to weapon target assignment and vehicle routing problems.

\begin{problem}[Generalized Assignment Problem, GAP] \label{prob:minSumForm}
    \begin{align*}
        \underset{\tilde{y}}{\text{minimize}} & \enskip c^{T}\tilde{y}\\
        \text{subject to} & \enskip \sum_{j=1}^{q}\tilde{y}_{ij}\leq1 \text{ for all }i=1,...,p\\
        & \enskip \sum_{i=1}^{p}t_{ij}\tilde{y}_{ij} \leq T_{\max, j}\text{ for all } j = 1,...,q\\
        & \enskip \tilde{y} \in \mathbb{Z}^{pq}
        \end{align*}
        where~$t_{ij} \in \mathbb{R}$ is a weight and~$T_{\max,j}$ is a~$p \times 1$ vector of each minimum value
        for resources per agent. 
\end{problem}

We reformulate the problem into a Kronecker Product form for condensed matrix notation.
\begin{problem}[GAP, Kronecker Product Form] \label{prob:kroneckerForm}
\begin{align}
\underset{\tilde{y}}{\text{minimize}} & \enskip c^{T}\tilde{y}\\
\text{subject to} & \enskip (\textbf{1}_{p}\otimes I_{q\times q})^{T}\tilde{y} \leq\textbf{1}_{q}\\
& \enskip Dx \leq T_{\max, p}\\
& \enskip \tilde{y} \in \mathbb{Z}^{pq}
\end{align}
where~$\textbf{1}_{p} \in \mathbb{R}^{p}$ is a $p \times 1$ vector of ones, $I_{q \times q}\in \mathbb{R}^{q \times q}$ is the~$q \times q$ identity matrix, 
$t_{ij}$ and $T_{\max, p}$ are defined as in Problem~\ref{prob:minSumForm}, 
and $D\in \mathbb{R}^{p\times pq}$ is a block diagonal matrix where each row is the price of the resource task from $\tilde{y}_{p1}, \tilde{y}_{p2}... \tilde{y}_{pq}$ for each agent denoted with the row vector $t_{p}$.
\end{problem}

Problem~\ref{prob:kroneckerForm} fits the form of Problem~\ref{prob:Original-MILP}, and we relax, regularize, and then solve it with Algorithm~\ref{alg:algorithm1}. 
In the purely integer case of Problem~\ref{prob:kroneckerForm}, we
have 
$\omega= \textbf{1}$, where $\textbf{1}$ is a vector of ones. Normally distributed random numbers are used to construct the cost vector $c$ for the generalized assignment problem.
Both varying communication rates and varying computation rates
were used in simulation with $100$ agents, $100$ tasks, $100$ primal agents, and $70$ dual agents. Algorithm~\ref{alg:algorithm1} was run for $10^5$ iterations in each case. 

\begin{figure}[ht] 
    \centering
    \includegraphics[width=3.2in,draft=false]{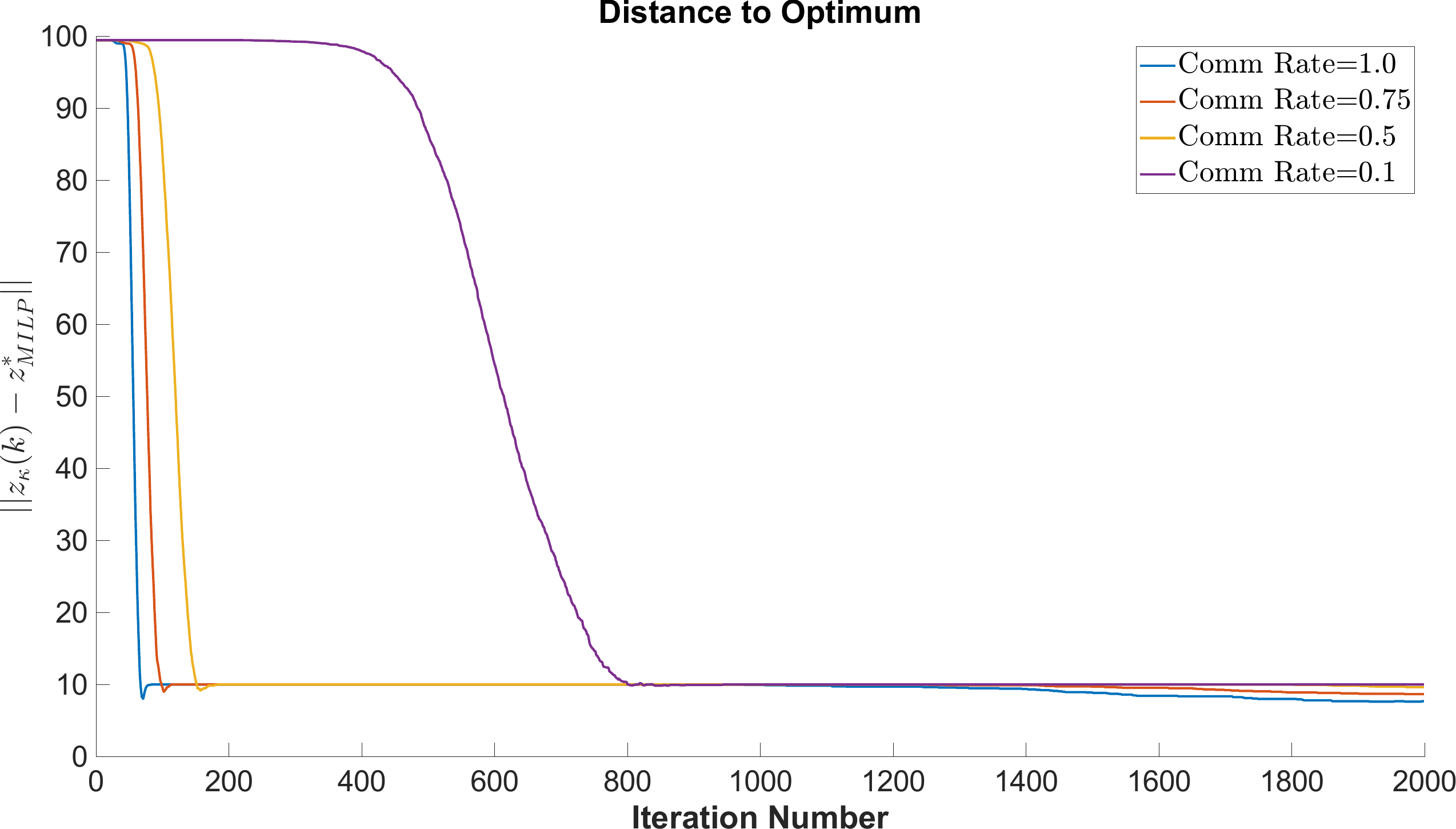}
    \caption{Four runs are compared based on varying communication rates for a random Generalized Assignment Problem with 100 agents and 100 tasks where there are 100 primal agents and 70 dual agents. The distance between iterates and the optimal solution obtained is compared with respect to the communication rates.
    As expected, a lower rate of communications results in agents requiring more iterations to converge. However, we do see that in all cases
    agents eventually reach a steady-state cost of approximately~$10$, followed by a period of slow decrease thereafter,
    indicating that agents reach the same solution, regardless of their communication rate. 
    }
    \label{fig:commRate}
    \end{figure}

From the simulations, it was observed that the iterates produced by Algorithm~\ref{alg:algorithm1} under all
communication rates produced iterates whose roundings 
satisfied all constraints after only $10^3$ iterations as seen in Figure~\ref{fig:commRate}. 

\begin{figure}[ht]
    \centering
    \includegraphics[width=3.25in,draft=false]{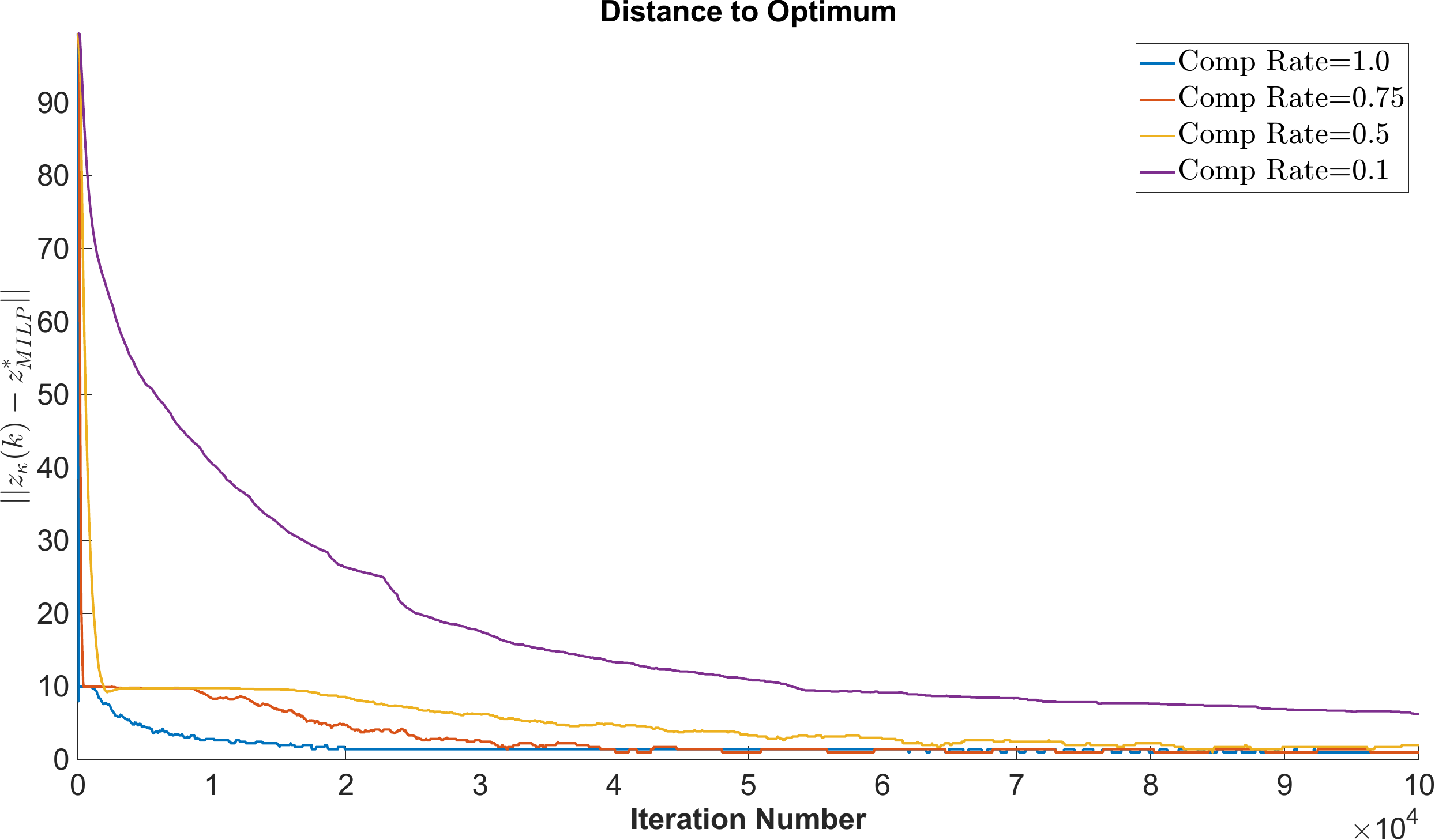}
    \caption{Four runs are compared based on varying computation rates for a random Generalized Assignment Problem with 100 agents and 100 tasks where there are 100 primal agents and 70 dual agents. The distance between iterates and the optimal solution obtained is compared with respect to the computation rates.
    As with changing communication rates, it is expected that agents with less frequent communications do indeed require
    more iterations to converge, which
    is what is seen here. 
    }
    \label{fig:compRate}
    \end{figure}

It was observed from Figure~\ref{fig:compRate} that only the computation rates of $1.0$ and $0.75$ 
produced iterates whose rounded forms satisfied all constraints after $10^5$ iterations. From Figure~\ref{fig:compRate}, 
computation rates of $0.5$ and $0.1$ both need longer than $10^5$ iterations 
to produce iterates whose roundings satisfy all constraints.
Nonetheless, they do indeed eventually produce iterates whose
rounded forms satisfy all constraints. 

\begin{figure}[ht] 
    \centering
    \includegraphics[width=3.25in,draft=false]{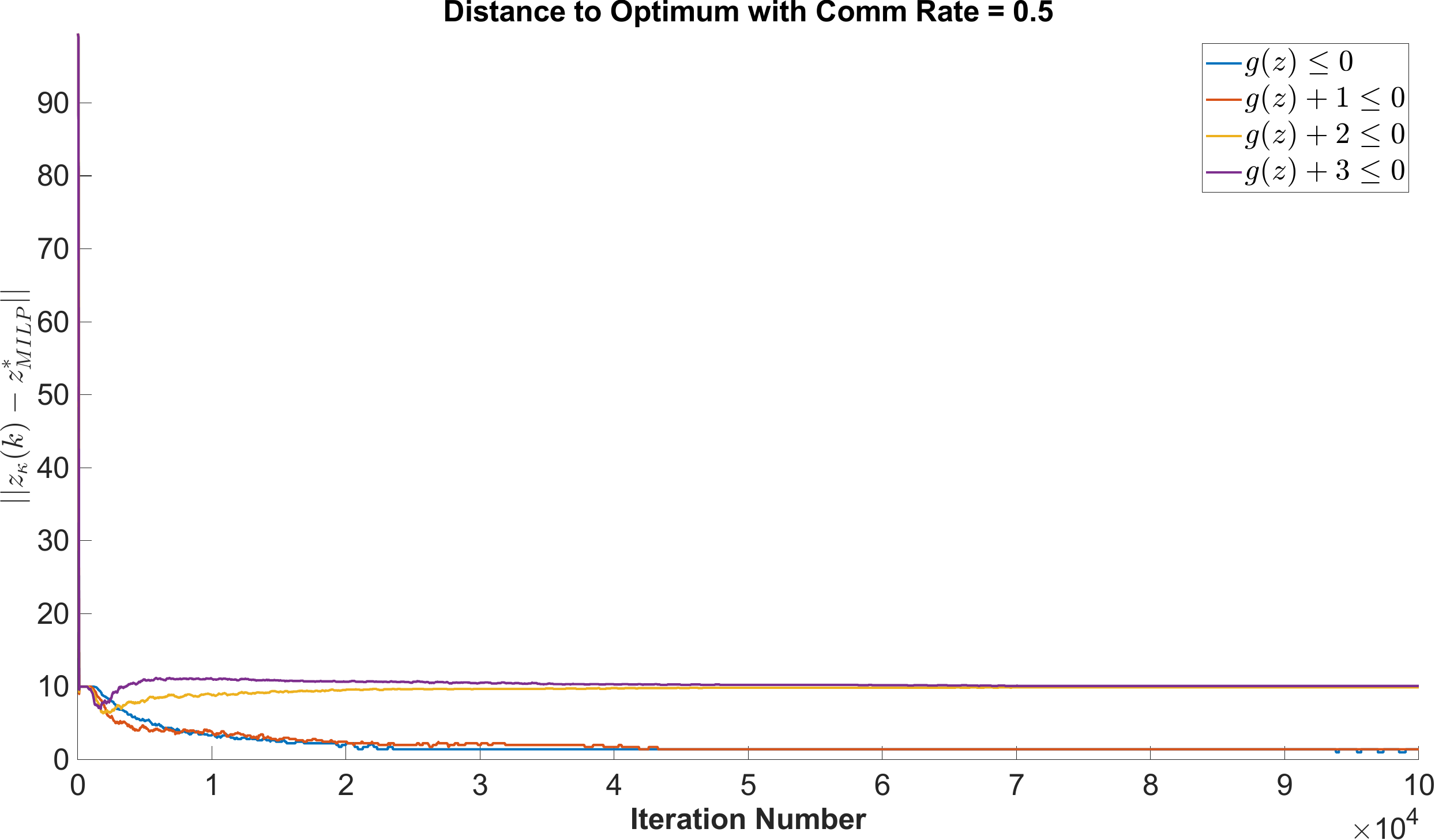}
    \caption{Four runs are compared based on varying 
    levels of constraint tightening. 
    The distance between iterates and the optimal solution obtained is compared with respect to the perturbation terms with a communication rate of $0.5$.
    As the constraints are tightened more, we see that agents have a higher cost at the final iterate. This is intuitive
    because tightening constraints leads to a smaller feasible region to optimize over, which necessarily leads
    to higher costs. 
    }
    \label{fig:slaterExample}
    \end{figure}

Additionally, Figure~\ref{fig:slaterExample} illustrates the impact of tightening constraints on the feasibility of the final output from the algorithm. As the term~$\phi$ increases, Algorithm~\ref{alg:algorithm1}'s final output
is further from the optimal solution. 

\section{Conclusion\label{sec:Conclusion}}
This paper relaxes and regularizes MILPs, then solves them in a distributed way and rounds the solution. 
This distributed algorithm was shown to be robust to asynchrony, and, 
despite both relaxation and regularization, feasibility of agents' final answer was guaranteed.
Moreover, theoretical bounds on
suboptimality in their costs were presented. 
Future work will address broader classes of mixed-integer programs and other classes of algorithms
for solving them.

\bibliographystyle{unsrt}
\bibliography{references}

\begin{thebibliography}{10}

\bibitem{richards2005mixed}
Arthur Richards and Jonathan How.
\newblock Mixed-integer programming for control.
\newblock {\em American Control Conference}, pages 2676--2683, 2005.

\bibitem{fischetti2018deep}
Matteo Fischetti and Jason Jo.
\newblock Deep neural networks and mixed integer linear optimization.
\newblock {\em Constraints}, 23(3):296--309, 2018.

\bibitem{hendricksondecentralized}
Katherine Hendrickson, Prashant Ganesh, Kyle Volle, Paul Buzaud, Kevin Brink,
  and Matthew Hale.
\newblock Decentralized weapon--target assignment under asynchronous
  communications.
\newblock {\em Journal of Guidance, Control, and Dynamics}, pages 1--13, 2022.

\bibitem{ioan2021mixed}
Daniel Ioan, Ionela Prodan, Sorin Olaru, Florin Stoican, and Silviu-Iulian
  Niculescu.
\newblock Mixed-integer programming in motion planning.
\newblock {\em Annual Reviews in Control}, 51:65--87, 2021.

\bibitem{conforti2014integer}
Michele Conforti, G{\'e}rard Cornu{\'e}jols, Giacomo Zambelli, et~al.
\newblock {\em Integer programming}, volume 271.
\newblock Springer, 2014.

\bibitem{neumann2019feasible}
Christoph Neumann, Oliver Stein, and Nathan Sudermann-Merx.
\newblock A feasible rounding approach for mixed-integer optimization problems.
\newblock {\em Computational Optimization and Applications}, 72:309--337, 2019.

\bibitem{camisa2021distributed}
Andrea Camisa, Ivano Notarnicola, and Giuseppe Notarstefano.
\newblock Distributed primal decomposition for large-scale milps.
\newblock {\em Transactions on Automatic Control}, 67(1):413--420, 2021.

\bibitem{camisa2019primal}
Andrea Camisa and Giuseppe Notarstefano.
\newblock Primal decomposition and constraint generation for asynchronous
  distributed mixed-integer linear programming.
\newblock {\em European Control Conference}, pages 77--82, 2019.

\bibitem{hendrickson2021totally}
Katherine~R Hendrickson and Matthew~T Hale.
\newblock Totally asynchronous primal-dual convex optimization in blocks.
\newblock {\em Transactions on Control of Network Systems}, 10(1):454--466,
  2022.

\bibitem{hale2017asynchronous}
Matthew~T Hale, Angelia Nedi{\'c}, and Magnus Egerstedt.
\newblock Asynchronous multiagent primal-dual optimization.
\newblock {\em Transactions on Automatic Control}, 62(9):4421--4435, 2017.

\bibitem{falsone2017dual}
Alessandro Falsone, Kostas Margellos, Simone Garatti, and Maria Prandini.
\newblock Dual decomposition for multi-agent distributed optimization with
  coupling constraints.
\newblock {\em Automatica}, 84:149--158, 2017.

\bibitem{vujanic2016decomposition}
Robin Vujanic, Peyman~Mohajerin Esfahani, Paul~J Goulart, S{\'e}bastien
  Mari{\'e}thoz, and Manfred Morari.
\newblock A decomposition method for large scale milps, with performance
  guarantees and a power system application.
\newblock {\em Automatica}, 67:144--156, 2016.

\bibitem{falsone2020distributed}
Alessandro Falsone and Maria Prandini.
\newblock A distributed dual proximal minimization algorithm for
  constraint-coupled optimization problems.
\newblock {\em Control Systems Letters}, 5(1):259--264, 2020.

\bibitem{camisa2019primalasync}
Andrea Camisa and Giuseppe Notarstefano.
\newblock Primal decomposition and constraint generation for asynchronous
  distributed mixed-integer linear programming.
\newblock {\em European Control Conference}, pages 77--82, 2019.

\bibitem{bragin2020distributedasync}
Mikhail~A Bragin, Bing Yan, and Peter~B Luh.
\newblock Distributed and asynchronous coordination of a mixed-integer linear
  system via surrogate lagrangian relaxation.
\newblock {\em Transactions on Automation Science and Engineering},
  18(3):1191--1205, 2020.

\bibitem{wolsey2020integer}
Laurence~A Wolsey.
\newblock {\em Integer programming}.
\newblock John Wiley \& Sons, 2020.

\bibitem{glover1975pitfalls}
Fred Glover and David~C Sommer.
\newblock Pitfalls of rounding in discrete management decision problems.
\newblock {\em Decision Sciences}, 6(2):211--220, 1975.

\bibitem{bertsimas2005optimization}
Dimitris Bertsimas and Robert Weismantel.
\newblock {\em Optimization over integers}, volume~13.
\newblock Dynamic Ideas, 2005.

\bibitem{bertsekas2003convex}
Dimitri Bertsekas, Angelia Nedic, and Asuman Ozdaglar.
\newblock {\em Convex analysis and optimization}, volume~1.
\newblock Athena Scientific, 2003.

\bibitem{arrow1958studies}
Kenneth~Joseph Arrow, Leonid Hurwicz, Hirofumi Uzawa, Hollis~Burnley Chenery,
  Selmer Johnson, and Samuel Karlin.
\newblock {\em Studies in linear and non-linear programming}, volume~2.
\newblock Stanford University Press Stanford, 1958.

\bibitem{aze2002sensitivity}
Dominique Az{\'e} and Jean-No{\"e}l Corvellec.
\newblock On the sensitivity analysis of hoffman constants for systems of
  linear inequalities.
\newblock {\em SIAM Journal on Optimization}, 12(4):913--927, 2002.

\bibitem{stein2016error}
Oliver Stein.
\newblock Error bounds for mixed integer linear optimization problems.
\newblock {\em Mathematical Programming}, 156:101--123, 2016.

\bibitem{bertsekas2015parallel}
Dimitri Bertsekas and John Tsitsiklis.
\newblock {\em Parallel and distributed computation: numerical methods}.
\newblock Athena Scientific, 2015.

\bibitem{ubl19}
Matthew Ubl and Matthew~T. Hale.
\newblock Totally asynchronous distributed quadratic programming with
  independent stepsizes and regularizations.
\newblock In {\em 2019 IEEE 58th Conference on Decision and Control (CDC)},
  pages 7423--7428, 2019.

\bibitem{hochhaus18}
Stefan Hochhaus and Matthew~T. Hale.
\newblock Asynchronous distributed optimization with heterogeneous
  regularizations and normalizations.
\newblock In {\em 2018 IEEE Conference on Decision and Control (CDC)}, pages
  4232--4237, 2018.

\bibitem{tseng1991rate}
Paul Tseng.
\newblock On the rate of convergence of a partially asynchronous gradient
  projection algorithm.
\newblock {\em SIAM Journal on Optimization}, 1(4):603--619, 1991.

\bibitem{koshal2011multiuser}
Jayash Koshal, Angelia Nedi{\'c}, and Uday~V Shanbhag.
\newblock Multiuser optimization: Distributed algorithms and error analysis.
\newblock {\em SIAM Journal on Optimization}, 21(3):1046--1081, 2011.

\end{thebibliography}

    \newpage
    \appendix

    \subsection{Proof of Theorem~\ref{thm:SlaterConnectionToGranularity}} \label{ss:theorem1}
    From \cite[Proposition 2.9]{neumann2019feasible}, the MILP 
    in Problem~\ref{prob:Original-MILP} 
    is granular if and only if the optimal cost in the following feasibility problem is less than one: 
    \begin{align*}
                \textnormal{Problem FP: }\underset{(x, y, \xi) \in X \times \mathbb{R}^{m} \times \mathbb{R}}{\text{minimize}} 
                & \enskip \xi \\
                \text{subject to} 
                 \enskip Ex + &Fy - \lfloor h \rfloor_{\omega} - \xi\omega + \frac{1}{2}\rho \leq 0 \\
                & y - y^u - \left(\xi - \frac{1}{2}\right)\bone \leq 0 \\
                & {-y} + y^{\ell} + \left(\frac{1}{2} - \xi\right)\bone \leq  0 \\
                &\xi \geq 0, 
                \end{align*}
                where all terms
                in Problem FP are defined
                in Problem~\ref{prob:matrix_Feasible_LP}.
                If we choose a constant positive value for~$\xi$ that is less than~$1$, then we can solve this problem by optimizing only over~$(x,y)$. 
                Using~$\xi_e$ from Proposition~\ref{prop:GranRounding}, we
                have~$\xi_e < 1$, and we
                consider~$\xi^* = \frac{1}{2}(\xi_e + 1)$,
                which satisfies~$\xi^* \in (\xi_e, 1)$.

                We can therefore solve the problem
                \begin{align*}
                    \textnormal{Problem FP$^*$: } \underset{(x,y) \in X \times \mathbb{R}^{m}}{\text{minimize}} 
                    & \enskip \xi^*\\
                    \text{subject to} 
                    \enskip Ex + &Fy - \lfloor h \rfloor_{\omega} - \xi\omega + \frac{1}{2}\rho \leq 0 \\
                    & y - y^u - \left(\xi^* - \frac{1}{2}\right)\bone \leq 0 \\
                    & {-y} + y^{\ell} + \left(\frac{1}{2} - \xi^*\right)\bone \leq 0. 
                    \end{align*}
    This problem clearly has optimal cost less than~$1$, and thus it
    is solved by any point~$(x,y)$ that satisfies the constraints. 
    By Remark~\ref{rem:slater}, the
    satisfaction of Assumption~\ref{as:slater} by~$g_{\xi_e}$ implies
    that Slater's condition also holds for~$g_{\xi^*}$, i.e.,
    the point~$(\bar{x}, \bar{y}) \in X \times \mathbb{R}^m$
    from Assumption~\ref{as:slater} also
    satisfies~$g_{\xi^*}(\bar{x}, \bar{y}) < 0$. 
    Then the feasible region of Problem FP$^*$ is non-empty,
    and in particular a solution exists to Problem FP$^*$. 
    Then, by~\cite[Proposition 2.9]{neumann2019feasible}
    we conclude that Problem~\ref{prob:Original-MILP} 
    is granular. 
    
    Next,
    Assumption~\ref{as:Xcompact} ensures compactness of~$X$,
    and, for any~$\xi \in [\xi_e, 1)$, Assumption~\ref{as:slater} ensures
    that the set
    \begin{equation}
        Y_{\xi} = \left\{y \in \mathbb{R}^m \mid 
        y^{\ell} + \left(\frac{1}{2} - \xi\right)\bone \leq y \leq 
    y^u + \left(\xi - \frac{1}{2}\right)\bone
        \right\}
    \end{equation}
    is non-empty; it is compact by inspection. 
    Then Problem~\ref{prob:matrix_Feasible_LP}
    is equivalent to an inequality constrained linear program
    over a compact set~$X \times Y_{\xi}$, and
    it has a non-empty feasible region. 
    Then it has at least one solution~$(x^*, y^*)$. 
    Finally, because Problem~\ref{prob:matrix_Feasible_LP}
    uses~$\xi \in [\xi_e, 1)$, Proposition~\ref{prop:GranRounding} 
    guarantees that~$(x^*, \round{y^*})$ 
    is feasible for Problem~\ref{prob:Original-MILP}.

 \subsection{Proof of Lemma~\ref{thm:newProblem}} \label{app:tightening}
    Let~$\hat{z}_{\alpha,\delta}$ denote the saddle point
        of~$L_{\alpha,\delta}$; it is unique
        because~$L_{\alpha,\delta}$ is strongly convex
        in~$z$ and strongly concave in~$\lambda$. 
        From \cite[Lemma 3.3]{koshal2011multiuser}, the maximum violation of the~$\ell^{th}$ constraint 
        by~$\hat{z}_{\alpha,\delta}$ is
        \begin{equation}
            \max\big\{0,g_{\ell}(\hat{z}_{\alpha,\delta})\big\} \leq \underset{z \in Z}{\max}\vert\vert\nabla g_{\ell}(z)\vert \vert_{2} \cdot \underset{\lambda \in \Lambda}{\max}\vert\vert \lambda\vert \vert_{2} \sqrt{\frac{\delta}{2\alpha}}. 
        \end{equation}
        By definition, 
            $\underset{z \in Z}{\max}\vert\vert\nabla g_{\ell}(z)\vert \vert_{2}  = \|A_{\ell}\|_2$.
        Using the fact that~$\|\lambda\|_2 \leq \|\lambda\|_1$, we can apply   
        Lemma~\ref{lem:saddlePoint} to upper bound~$\|\lambda\|_2$,
        which completes the proof. 
        %

    \subsection{Proof of Lemma~\ref{lem:saddlePoint}}\label{ss:lem1}
Using inequality in~\eqref{eq:saddleineqs}, we see that
for~$\bar{z}$ we have
\begin{equation}
L_{\kappa}(\hat{z}_{\kappa}, 0) \leq L_{\kappa}(\hat{z}_{\kappa}, \hat{\lambda}_{\kappa}) \leq
L_{\kappa}(\bar{z}, \hat{\lambda}_{\kappa}),
\end{equation}
where expanding the definition of~$L_{\kappa}$ 
from~\eqref{eq:RegularizedLagran} gives
\begin{equation}
f(\hat{z}_{\kappa}) + \frac{\alpha}{2}\|\hat{z}_{\kappa}\|^2 \leq 
f(\bar{z}) + \frac{\alpha}{2}\|\bar{z}\|^2 + 
\hat{\lambda}_{\kappa}\big(A\bar{z} - b - \nu -\phi\bone) - \frac{\delta}{2}\|\hat{\lambda}_{\kappa}\|^2. 
\end{equation}
Since~$\alpha > 0$ and~$\delta > 0$ we then have
\begin{equation}
f(\hat{z}_{\kappa})  \leq 
f(\bar{z}) + \frac{\alpha}{2}\|\bar{z}\|^2 + 
\hat{\lambda}_{\kappa}\big(A\bar{z} - b - \nu -\phi\bone). 
\end{equation}
With~$B_0(r)$ as defined in~\eqref{eq:b0}, we have
$\min_{z \in B_0(r)} f(z) \leq f(\hat{x}_{\kappa})$ and thus
\begin{equation}
\min_{z \in B_0(r)} f(z) \leq f(\bar{z}) + \frac{\alpha}{2}\|\bar{z}\|^2 + 
\hat{\lambda}_{\kappa}\big(A\bar{z} - b - \nu-\phi\bone).
\end{equation}

We rearrange terms to find
\begin{align}
f(\bar{z}) + \frac{\alpha}{2}\|\bar{z}\|^2\! - \!\!\!\min_{z \in B_0(r)} \!f(z) 
&\geq
-\!\sum_{j=1}^{p} \hat{\lambda}_{\kappa,j} \big(A_j\bar{z} \!-\! b_j \!-\! \nu_j \!-\! \phi\big) \nonumber \\
&\geq {\min_{1 \leq \ell \leq p}} \big({-A}_{\ell}\bar{z} + b_{\ell} + \nu_{\ell} + \phi\big) \sum_{j=1}^{p} \hat{\lambda}_{\kappa,j}. 
\end{align}



We have~$-\big(A_{\ell}\bar{z} - b_{\ell} - \nu_{\ell}-\phi\big) > 0$
 because~$\bar{z}$ is a Slater point,
and we can divide by the~$\min$ term without reversing the inequality. 
Then 
\begin{equation}
\hat{\lambda}_{\kappa}\in \Lambda \!=\! \left\{\!\lambda \!\in\! \mathbb{R}_{+}^{y} \!:\! \vert\vert\lambda\vert\vert_{1} \!\leq\! \frac{f(\bar{z}) \!+\! \frac{\alpha}{2}\vert\vert\bar{z}\vert\vert^{2} -\!\!\!\! \underset{z\in B_0(r)}{\min}f(z)}{{\min\limits_{1\leq \ell \leq m}} 
{-A}_{\ell}\bar{z} + b_{\ell} + \nu_{\ell} + \phi}\right\}. \label{eq:Lem4_part1}
\end{equation}
Finally, we note that
\begin{equation}
\underset{z\in B_{o}(r)}{\min}c^{T}z=c^{T}\cdot\begin{pmatrix}-\frac{c}{\vert\vert c\vert\vert}\end{pmatrix}\cdot r=-\vert\vert c\vert\vert\cdot r.\label{eq:Lem4_part3}
\end{equation}
Substituting \eqref{eq:Lem4_part3}
into \eqref{eq:Lem4_part1} gives the desired bound.

\subsection{Proof of Lemma~\ref{lem:regError}} \label{app:regError}
Using \cite[Lemma 3.3]{koshal2011multiuser},
we can upper bound between the difference in costs at $\hat{z}_{\kappa}$ and
$z^{*}$ via
\begin{equation}
\vert f(\hat{z}_{\kappa})-f(z^{*})\vert\leq\underset{z\in Z}{\max}\vert\vert\nabla f(z)\vert\vert\underset{\lambda\in \Lambda}{\max}\vert\vert\lambda\vert\vert\sqrt{\frac{\delta}{2\alpha}}+\frac{\alpha}{2}\underset{z\in Z}{\max}\vert\vert z\vert\vert.\label{eq:optBoundFull}
\end{equation}
From Problem \ref{prob:matrix_Feasible_LP}, 
$\underset{z\in Z}{\max}\vert\vert\nabla f(z)\vert\vert=\vert\vert c\vert\vert$
since $\nabla f \equiv c.$ Next, from~\eqref{eq:r} we have
\begin{equation}\label{eq:maxR}
\underset{z\in Z}{\max}\vert\vert z\vert\vert = r. 
\end{equation}
Since $\ensuremath{\lambda}$ is a vector,
we always have $\underset{\lambda\in \Lambda}{\max}\vert\vert\lambda\vert\vert_{2}\leq\underset{\lambda\in \Lambda}{\max}\vert\vert\lambda\vert\vert_{1}.$
From Lemma \ref{lem:saddlePoint}, we can therefore upper bound $\vert\vert\lambda\vert\vert_{2}$
as
\begin{equation}
 \underset{\lambda\in \Lambda}{\max}\vert\vert\lambda\vert\vert_{2} \leq \frac{c^{T}\bar{z}+\frac{\alpha}{2}\vert\vert\bar{z}\vert\vert^{2} + \vert\vert c\vert\vert\cdot r}{\min\limits_{1\leq j\leq m} -A_{j}\bar{z} + b_{j} + \nu_{j} + \phi}.\label{eq:boundmaxLamb-1}
\end{equation}
Using \eqref{eq:maxR} and \eqref{eq:boundmaxLamb-1} in \eqref{eq:optBoundFull} completes the proof. 
\hfill $\blacksquare$

\subsection{Proof of Lemma~\ref{lem:(Dual-Convergence-Between}}\label{ss:lem7}
Define $\hat{z}_{\kappa}(t_{1}B) = \argmin_{z\in Z}L_{\kappa}(z,\lambda(t_{1}B))$
and recall $\hat{z}_{\kappa} = \argmin_{z\in Z}L_{\kappa}(z,\hat{\lambda}_{\kappa})$. 
Given that all dual agents use the same primal variables in
their computations, we analyze all dual agents' computations
simultaneously with the combined dual update law
\begin{equation}
\lambda(t_2B) = \Pi_{\Lambda}\big[\lambda(t_1B) + \beta \nabla_{\lambda}L_{\kappa}\big(z(t_2B), \lambda(t_1B)\big)\big],
\end{equation}
where~$z(t_2B)$ is the common primal variable that the dual agents use to compute~$\lambda(t_2B)$. 
Expanding the dual update law and using the non-expansiveness
of $\Pi_{\Lambda}$ gives
\begin{multline}
\vert\vert\lambda(t_{2}B)-\hat{\lambda}_{\kappa}\vert\vert  
=(1-\beta\delta)^{2}\vert\vert\lambda(t_{1}B)-\hat{\lambda}_{\kappa}\vert\vert^{2} 
 +\beta^{2}\vert\vert g(z(t_{2}B)-g(\hat{z}_{\kappa})\vert\vert^{2} \\ 
-2\beta(1-\beta\delta)(\lambda(t_{1}B)-\hat{\lambda}_{\kappa})^{T}(g(\hat{z}_{\kappa})-g(z(t_{2}B))).
\end{multline}

Adding $g(\hat{z}_{\kappa}(t_{1}B))-g(\hat{z}_{\kappa}(t_{1}B))$
inside the last set of parentheses gives
\begin{multline}
\vert\vert\lambda(t_{2}B)-\hat{\lambda}_{\kappa}\vert\vert = (1-\beta\delta)^{2}\vert\vert\lambda(t_{1}B)-\hat{\lambda}_{\kappa}\vert\vert^{2} 
 - 2\beta(1-\beta\delta)(\lambda(t_{1}B)-\hat{\lambda}_{\kappa})^{T}(g(\hat{z}_{\kappa})-g(\hat{z}_{\kappa}(t_{1}B)) \\
 - 2\beta(1-\beta\delta)(\lambda(t_{1}B)-\hat{\lambda}_{\kappa})^{T}(g(\hat{z}_{\kappa}(t_{1}B))-g(z(t_{2}B))) 
 +\beta^{2}\vert\vert g(z(t_{2}B)-g(\hat{z}_{\kappa})\vert\vert^{2} .\label{eq:Lam_t2_LamOpt}
\end{multline}

Applying~\cite[Lemma~4.1]{koshal2011multiuser} 
to the pairs $(\hat{z}_{\kappa},\hat{\lambda}_{\kappa})$
and $(\hat{z}_{\kappa}(t_{1}B),\lambda(t_{1}B))$ results in 
\begin{equation}
(\lambda(t_{1}B)-\hat{\lambda}_{\kappa})^{T}(-g(\hat{z}_{\kappa})+g(\hat{z}_{\kappa}(t_{1}B)))  
  \geq\frac{\alpha}{\vert\vert A\vert\vert_{2}^{2}}\vert\vert g(\hat{z}_{\kappa}(t_{1}B))-g(\hat{z}_{\kappa})\vert\vert^{2},\label{eq:applyingKoshLemma}
\end{equation}
which we apply to the third term on
the right-hand side of \eqref{eq:Lam_t2_LamOpt} to find 
\begin{multline}
\vert\vert\lambda(t_{2}B) - \hat{\lambda}_{\kappa}\vert\vert  \leq(1-\beta\delta)^{2}\vert\vert\lambda(t_{1}B)-\hat{\lambda}_{\kappa}\vert\vert^{2}
 - 2\beta(1-\beta\delta)\frac{\alpha}{\vert\vert A\vert\vert_{2}^{2}}\vert\vert g(\hat{z}_{\kappa}(t_{1}B))-g(\hat{z}_{\kappa})\vert\vert^{2} \\
 - 2\beta(1-\beta\delta)(\lambda(t_{1}B)-
 \hat{\lambda}_{\kappa})^{T}(g(\hat{z}_{\kappa}(t_{1}B))-g(z(t_{2}B))) 
  + \beta^{2}\vert\vert g(z(t_{2}B)-g(\hat{z}_{\kappa})\vert\vert^{2}. \label{eq:Lam_t2_LamOpt-1}
\end{multline}

To bound the last term, we see that $0\leq\vert\vert(1-\beta\delta)(g(\hat{z}_{\kappa}(t_{1}B))-g(z(t_{2}B)))+
\beta(\lambda(t_{1}B)-\hat{\lambda}_{\kappa})\vert\vert^{2}$, which
can be expanded and rearranged to give
\begin{multline}
-2\beta(1-\beta\delta)(\lambda(t_{1}B)-\hat{\lambda}_{\kappa})^{T}(g(\hat{z}_{\kappa}(t_{1}B))-g(z(t_{2}B))) \leq \\ (1-\beta\delta)^{2}\vert\vert(g(\hat{z}_{\kappa}(t_{1}B))-g(z(t_{2}B)))\vert\vert^{2} +\beta^{2}\vert\vert(\lambda(t_{1}B)-\hat{\lambda}_{\kappa})\vert\vert^{2}.\label{eq:lastTermBound}
\end{multline}

Then applying \eqref{eq:lastTermBound} to the
last term in \eqref{eq:Lam_t2_LamOpt-1} to obtain
\begin{multline}
\vert\vert\lambda(t_{2}B) - \hat{\lambda}_{\kappa}\vert\vert \leq(1-\beta\delta)^{2}\vert\vert\lambda(t_{1}B)-\hat{\lambda}_{\kappa}\vert\vert^{2} 
 +\beta^{2}\vert\vert g(z(t_{2}B)-g(\hat{z}_{\kappa})\vert\vert^{2} 
  -2\beta(1-\beta\delta)\frac{\alpha}{\vert\vert A\vert\vert_{2}^{2}}\vert\vert g(\hat{z}_{\kappa}(t_{1}B))-g(\hat{z}_{\kappa})\vert\vert^{2} 
 \\ +(1-\beta\delta)^{2}\vert\vert(g(\hat{z}_{\kappa}(t_{1}B))-g(z(t_{2}B)))\vert\vert^{2} 
 +\beta^{2}\vert\vert(\lambda(t_{1}B)-\hat{\lambda}_{\kappa})\vert\vert^{2}.\label{eq:Lam_t2_LamOpt-FinalBound}
\end{multline}

Next, add and subtract $g(\hat{z}_{\kappa}(t_{1}B))$ inside the norm
in the second term of \eqref{eq:Lam_t2_LamOpt-FinalBound} to obtain
\begin{multline}
\beta^{2}\vert\vert g(z(t_{2}B))-g(\hat{z}_{\kappa})\vert\vert^{2} \leq\beta^{2}\vert\vert g(z(t_{2}B))-g(\hat{z}_{\kappa}(t_{1}B))\vert\vert^{2}  +\beta^{2}\vert\vert g(\hat{z}_{\kappa}(t_{1}B))-g(\hat{z}_{\kappa})\vert\vert^{2}\\
 + 2\beta^{2}\vert\vert g(z(t_{2}B))-g(\hat{z}_{\kappa}(t_{1}B))\vert\vert\cdot\vert\vert g(\hat{z}_{\kappa}(t_{1}B))-g(\hat{z}_{\kappa})\vert\vert.
\end{multline}

Applying this to the second term in \eqref{eq:Lam_t2_LamOpt-FinalBound}
gives
\begin{multline}
\vert\vert\lambda(t_{2}B) - \hat{\lambda}_{\kappa}\vert\vert \leq((1-\beta\delta)^{2}+\beta^{2})\vert\vert\lambda(t_{1}B)-\hat{\lambda}_{\kappa}\vert\vert^{2} 
+\begin{pmatrix}\beta^{2}-2\beta(1-\beta\delta)\frac{\alpha}{\vert\vert A\vert\vert_{2}^{2}}\end{pmatrix}\vert\vert g(\hat{z}_{\kappa}(t_{1}B))-g(\hat{z}_{\kappa})\vert\vert^{2}\\
 +((1-\beta\delta)^{2}+\beta^{2})\vert\vert g(\hat{z}_{\kappa}(t_{1}B))-g(z(t_{2}B))\vert\vert^{2} 
 +2\beta^{2}\vert\vert g(z(t_{2}B)-g(\hat{z}_{\kappa}(t_{1}B))\vert\vert\vert\vert g(\hat{z}_{\kappa}(t_{1}B))-g(\hat{z}_{\kappa})\vert\vert.
\end{multline}

By hypothesis $0<\beta<\frac{2\alpha}{\vert\vert A\vert\vert_{2}+2\alpha\delta},$ which
makes the second term negative. Dropping this negative term and applying
the Lipschitz property of $g$ give the upper bound
\begin{multline}
\vert\vert\lambda(t_{2}B)-\hat{\lambda}_{\kappa}\vert\vert \leq((1-\beta\delta)^{2}+\beta^{2})\vert\vert\lambda(t_{1}B)-\hat{\lambda}_{\kappa}\vert\vert^{2} 
+((1-\beta\delta)^{2}+\beta^{2})\vert\vert A\vert\vert_{2}^{2}\vert\vert\hat{z}_{\kappa}(t_{1}B)-z(t_{2}B)\vert\vert^{2}\\
 +2\beta^{2}\vert\vert A\vert\vert_{2}^{2}\vert\vert z(t_{2}B)-\hat{z}_{\kappa}(t_{1}B)\vert\vert\cdot\vert\vert\hat{z}_{\kappa}(t_{1}B)-\hat{z}_{\kappa}\vert\vert. 
\end{multline}

Lemma \ref{lem:primalConvergence} can be applied to $\vert\vert\hat{z}_{\kappa}(t_{1}B)-z(t_{2}B)\vert\vert^{2}$
to give an upper bound of
\[
\vert\vert\hat{z}_{\kappa}(t_{1}B)-z(t_{2}B)\vert\vert^{2}\leq q_{p}^{2(t_{2}-t_{1})}\vert\vert z(t_{1}B)-\hat{z}_{\kappa}(t_{1}B)\vert\vert^{2},
\]
where $q_{p}\text{ }=(1-\theta\gamma)\in[0,1)$. Next, the maximum
distance between any two primal variables is bounded by 
\begin{equation}
\underset{z,y\in Z}{\max}\vert\vert z-y\vert\vert \leq \max_{z, y} \|z\| + \|y\| \leq 2r.
\end{equation}
Using this maximum distance and setting $q_{d} = (1-\beta\delta)^{2}+\beta^{2}$,
it follows that $q_{d}\in[0,1)$ because $\beta<\frac{2\delta}{1+\delta^{2}}$, which
completes the proof.

\end{document}